\documentclass[reqno]{amsart}
\usepackage{amscd,amssymb,amsmath,amsfonts,amsthm,latexsym}
\usepackage{extpfeil}
\usepackage[T1]{fontenc}
\usepackage[all]{xy}

\usepackage{lscape}
\usepackage{graphicx}
\usepackage{xcolor}
\usepackage{cite}
\usepackage{multicol}
\usepackage{mathtools}
\usepackage{a4wide,hyperref}
\usepackage{enumitem}

\usepackage{listings}

\newdir{>}{!/8pt/\dir{}*\dir{>}}

\DeclareMathOperator{\lcm}{lcm}

\DeclareMathOperator{\Mat}{Mat}

\newtheorem{theorem}{Theorem}[section]
\newtheorem{proposition}[theorem]{Proposition}
\newtheorem{corollary}[theorem]{Corollary}
\newtheorem{lemma}[theorem]{Lemma}
\newtheorem*{lemm}{Lemma}
\newtheorem{definition}[theorem]{Definition}
\newtheorem{example}[theorem]{Example}

\newtheorem{remark}[theorem]{Remark}
\numberwithin{equation}{section}

\definecolor{Orange}{rgb}{1,0.5,0}

\begin{document}

\title{Gr\"obner Bases under Composition, Survey paper}
\author{Mahmoud S. Alsersawi}
\address{[Mahmoud S. Alsersawi]  University of Santiago de Compostela, 15782, Spain.
}
\email{aboyahya1986@hotmail.com}

\author{Manuel Ladra}
\address{[Manuel Ladra] Department of Mathem\'aticas, Institute of Matem\'aticas, University of Santiago de Compostela, 15782, Spain.}
\email{manuel.ladra@usc.es}

\thanks{Most of the results presented in this work were obtained thanks to the financing of the Programme \textbf{Erasmus Mundus Programme}:
Program for Excellence Academy Cooperation Exchange (PEACE II)-GA 2013--2443.}

\begin{abstract}
This paper is a survey on the study of the behaviour of the composition of polynomials on the computation of Gr\"obner bases. This survey brings together  some works published between 1995 and 2007. The authors of these papers gave answers to some questions in this subject for several types of Gr\"obner bases, over different monomials orderings and over different polynomials rings. Some of these answers are complete and some are not. Some papers of them were given to answer some open questions given in the first published paper of these papers and other some were given to generalize previous work. Note that, there are many other works in this subject over other special cases and types of Gr\"obner bases under the usual definition over arbitrary fields, over other types of bases, under the non-commutative case and with other methods of proofs and calculations.
\end{abstract}

\maketitle

\section{Introduction}
As many mathematicians know, the remainder of a polynomial $f$ on division by a set of finite polynomials $F$ using the Division Algorithm on the polynomial ring $K[x_1,\dots, x_n]$, where $K$ is a field and using some fixed monomial ordering depends on how we order the polynomials inside $F$ and thus this remainder is not unique. Many problems in theoretical and applied algebra obtained because of this reason, especially when we work with applications dealing with polynomial ideals, like the ideal description and membership problems and the problem of solving polynomial equations.

In 1965, using the Hilbert Basis Theorem, Bruno Buchberger gave in his PhD thesis \cite{BB} an algorithm to compute a basis for polynomial ideals called Gr\"obner basis. This basis depends in its definition on the monomial ideals, that is, it is a generating set of the ideal such that the set of its leading monomials generates the ideal generating from the set of the leading monomials of that ideal. This revolutionary idea put us on the road to solve many problems like what we considered before, because the remainder on division by a set which is a Gr\"obner basis will be unique now and does not depend on how we order of the polynomials inside this set. Gr\"obner basis after this became a very important tool in many branches of mathematics like computer algebra, geometry, etc. and have many applications in many fields like Elimination Theory. Many authors, including Buchberger himself, later started from this point and gave many modifications in the theory of Gr\"obner bases by giving several types of Gr\"obner bases with additional conditions like minimal, reduced and homogeneous Gr\"obner bases or by modifying algorithms to compute these types of Gr\"obner bases.

In 1995, the  Korean mathematician Hoon Hong started a subject on this theory when he asked the following question: what is the  behaviour of Gr\"obner bases with respect to some monomial ordering $<$ under the operation of polynomial composition. In other words, let $G$ be a Gr\"obner basis with respect to some monomial ordering $<$ for some ideal $I$, and let $\Theta$ be a list of $n$ polynomials on the polynomial ring $K[x_1,\dots, x_n]$, what are the conditions that must be put on the list $\Theta$ to make the set $G\circ\Theta$ be a Gr\"obner basis with respect to the same $<$ for the ideal $\langle I \circ \Theta \rangle$. He asked this question in the form: when does Gr\"obner basis computation commute with composition? He studied this problem and published a paper under the name ``Groebner Basis under Composition I''  \cite{Ho1}, and proved that the computation of a Gr\"obner bases with respect to some $<$  commutes with composition if and only if the composition is compatible with the monomial ordering $<$ and the non-divisibility. This paper was the first paper in this subject but not the last one, because Hong put on its end three open questions related to his work. Namely:

\begin{itemize}
	\item (Q1) Does there exist a decision procedure that will determine whether a given composition is compatible with a given term ordering. If so, find one.
	
	\item (Q2) When does a composition commute with the reduced Gr\"obner basis computation?
	
	\item (Q3) Let $G$ be a Gr\"obner basis for $I$ with respect to $<$. When $G \circ \Theta$ is a Gr\"obner basis for $\langle I \circ \Theta \rangle$  (possibly with respect to another monomial  ordering $<'$)?
\end{itemize}

These questions solved completely in the period 1996--2003 as we will explain later. Hong himself started this process in his second paper  \cite{Ho2} on this subject published in 1996, where he studied a special case with a special understanding of the words `another monomial ordering' in the third question where the other one may be different from the first one, but  he has a special way to define it using the first one itself and the list $\Theta$. In this paper, he studied this type of understanding with a special way of defining the other ordering using both of the first fixed one and the list $\Theta$ and gave the needed conditions to make this new ordering be a monomial ordering. He proved that if the list of leading monomials with respect to some $<$ of the list $\Theta$ is permuted powering, then for any Gr\"obner basis $G$ of an ideal $I$ with respect to the obtained ordering (which will be a monomial ordering under these conditions), the set  $G \circ \Theta$ is a Gr\"obner basis for $\langle I \circ \Theta \rangle$  with respect to $<$.

In 1998, the two Spanish mathematicians J. Guti\'errez and R. Rubio San Miguel gave a full answer for the second question about the behaviour of the composition on the computation of reduced Gr\"obner bases in  \cite{GR}. They proved that for every reduced Gr\"obner basis $G$ with respect to $<$, $G \circ\Theta$ is a reduced Gr\"obner basis with respect to the same monomial ordering $<$ if and only if the composition by $\Theta$ is compatible with the monomial ordering $<$ and $\Theta$ is a list of permuted univariate and monic polynomials. This result is different of Hong's result by the additional two conditions: univariate and monic. Every new condition will be needed for a reason comes from the meaning of the word `Reduced'. They also studied the special case of the third open question for reduced Gr\"obner bases and gave a sufficient condition to determine when the composition commutes with reduced Gr\"obner bases computation under some possibly different monomial orderings where the second possibly different monomial ordering is which the ordering defined before by Hong in \cite{Ho2}. In addition to these studies for reduced Gr\"obner bases case, they studied these two questions for minimal Gr\"obner bases case and proved some results close to the results proved by Hong in his two papers \cite{Ho1,Ho2}.

During the period where some authors tried to solve the open questions given by Hong, the Sweden mathematician P. Nordbeck studied this subject, but in some different way. Firstly, in 2001 in  \cite{Nor1}, where he studied the case of non-commutative Gr\"obner bases under composition. This is more difficult than the commutative case, but he used some of Hong's ideas and proved that the composition by $\Theta$ commutes with the computation of non-commutative Gr\"obner bases with respect to some $<$ if and only if the composition is compatible with the ordering $<$ and the set of leading words of $\Theta$ is combinatorially free. The second paper published \cite{Nor2} by him was in 2002, and he studied  the same question for some type of bases needed for studying the sub-algebras of polynomial rings and related to non-commutative polynomial rings called SAGBI bases. He gave a complete answer for this, and also  showed that this commutativity happened if and only if the composition is compatible with the fixed ordering.

After 2001, many Chinese mathematicians worked in this subject and published more than 9 papers (what are related directly to this subject) during the period 2001--2011. The most interesting thing is that J. Liu is a common name of the authors in 8 of these 9 papers. The only one without him was the paper  \cite{LZW} published by Z. Liu and M. Wang in 2001. In this paper, they solved the third question given by Hong in the general meaning of words `another monomial ordering',  which is when both of the two monomial orderings are arbitrary without any special conditions on any of them. Note also this way of generalization the previous results, but they worked under the general case where there are two possibly different polynomial rings, one before the composition and the other one after it where these ways of generalizations were studied together in this paper, that is the generalized definitions and results given in
it is by using  different polynomial rings $K[x_1, \dots, x_n]$ and $K[y_1, \dots, y_m]$ and different monomial orderings $<_1$ and $<_2$ before and after the composition.
They generalised all the results given by Hong in their result: the composition by $\Theta$ commutes with the computation of Gr\"obner bases with respect to $<_1$ and $<_2$ if and only if this composition  is compatible with $<_1$ and $<_2$ and for all $i\neq j$, the monomials $LM_{<_2}(\theta_i)$ and $LM_{<_2} (\theta_j)$ are relatively prime.
Also, they generalised the results given by J. Guti\'errez and R. Rubio San Miguel  for the reduced Gr\"obner bases case, and gave some good answers using different sufficient and necessary condition for this case of commutativity as we will see in Theorems~\ref{reduced Gr\"obner Bases Under Composition,diff monomial ordering, main result}, \ref{reduced Gr\"obner Bases Under Composition,diff monomial ordering, main result2}, and Corollary~\ref{reduced Gr\"obner Bases Under Composition,diff monomial ordering, coro}. This finished completely Hong's second and third questions.

The first paper of J. Liu in this subject was  \cite{LF} with X. Fu,  and published in 2002. In this paper, they gave a better answer for the Hong's first question asking about a decision procedure that will determine whether a given composition is compatible with a given monomial ordering or not. Their work used the concepts of linear algebra depending on the fact that every monomial ordering has a corresponding matrix proved with an algorithm of computing in many previous works like \cite{BeWe,DMY,Rob}. In 2003, J. Liu join with Z. Liu and M. Wang   published the paper   \cite{LZM}, where they completed and finished the work on Hong's first question using elementary rational row operations for matrices to obtain a decision procedure for this question.

After finishing these three open questions, some authors searched for other questions related to our subject. They used the same methods given by Hong to define and study other types of commutativity of the composition with the computation of other several types of Gr\"obner bases or other types of bases which are not necessarily Gr\"obner bases. The reader can notice the similarity of their definitions and main results and even the method of proofs and all these things are similar to those given by Hong. Some of these results will be like a generalization of Hong's results or some other previous results.

 During our research we studied  the following papers ordered historically: J. Liu, J.  Hou, Z. Peng and  W. Cao \cite{LHPC} in 2005
proved that the composition commutes with the computation of universal Gr\"obner bases (a universal Gr\"obner basis of the ideal $I$ is a Gr\"obner basis for $I$ for all monomial orderings) if and only if the composition is a single variable, i.e.  every single $\theta_i$ is a function in one variable that is not used in any other $\theta_j$.
J. Liu and M. Wang  \cite{LW1} studied in 2006 the question: what conditions we must put on the list $\Theta$ of homogeneous polynomials of the same degree to make every homogeneous Gr\"obner basis $G$ with respect to some monomial ordering $<$ for some homogeneous ideal $I$ compose using $\Theta$ to a homogeneous Gr\"obner basis $G\circ\Theta$ with respect to the same $<$ for the homogeneous ideal $ \langle I \circ \Theta \rangle$. They gave a complete answer: for every homogeneous Gr\"obner basis $G$ with respect to $<$, $G \circ\Theta$ is a homogeneous Gr\"obner basis with respect to the same monomial ordering if and only if the composition by $\Theta$ is homogeneously compatible with the monomial ordering $<$ and $\Theta$  is permuted powering.
	J. Liu and M. Wang in 2007 decided to generalise their results and the results proved by Hong in his first paper and unify these results in a common result in \cite{LW2}, where they studied the question of when homogeneous Gr\"obner bases under an arbitrary grading remain Gr\"obner bases after composition with respect to any fixed monomial ordering and gave a complete answer for it. In more details, let $\varGamma$ be an arbitrary grading on $K[x_1, \dots, x_n]$. If $G$ is any $\varGamma$-homogeneous Gr\"obner basis with respect to some monomial ordering $<$ for some ideal $I$, what are the conditions must the list $\Theta$ have to make the set $G\circ\Theta$ a Gr\"obner basis with respect to the same $<$ for the ideal $ \langle I \circ \Theta \rangle$. They gave a necessary and sufficient condition to the above problem, which is that the composition by $\Theta $ is $\varGamma$-compatible with the monomial ordering $<$ and the list $LM_{<}(\Theta)$ is permuted powering. If we choose some special values of $\varGamma$, then we will get the results proved in Hong's \cite{Ho1,LW1} and as we said above which may be regarded as a common generalization of these results. In this paper, no other special conditions will put on $\Theta$ since they studied the general case and the obtained composed set $G\circ\Theta$ needed to be only a Gr\"obner basis without any special cases.

In the papers \cite{LLC,TZ,ZAI,LLL1,LLL2}  the same subject was studied but over several special types of Gr\"obner bases or other types of bases using the same methods considered before and their main results are in the same way.

 In this paper we ignored many other papers in this subject and described a survey of some works published between 1995 and 2007. These papers have a common subject, which is to study the behaviour of the composition of polynomials on the computation of Gr\"obner bases. The authors of these papers gave answers to some questions in this subject for several types of Gr\"obner bases, over different monomials orderings and over different polynomials rings. Some of these answers are complete and some are not. Some papers of them were given to answer some open questions given in the first paper of them and other some were given to generalize previous work.
	
This survey  is contained in the second chapter of Mahmoud Alsersawi's doctoral thesis ``Gr\"obner Bases under Composition  over Fields with Valuations''  and is organized as follows.
In Section~\ref{S:Basic} we will give the most important definitions and results related to the basic theory of Gr\"obner bases.
	In the next two sections, Sections~\ref{Gr\"obner Bases Under Composition1} and \ref{Gr\"obner Bases Under Composition2}, we will study the work of Hong in his first paper  \cite{Ho1}, which was the first time to study the behaviour of the composition of polynomials on the computation of Gr\"obner bases. This paper appeared firstly in 1995 and another copy of it published in 1998. His way of definition and method of proof in this paper will be after it the guide of other authors who studied this subject later. At the end of this paper, Hong put three questions to be answered. Some answers will be studied in some following sections.
	
In Section~\ref{Gr\"obner Bases Under Composition3} we look into Hong's second paper \cite{Ho2} published in 1996, where  he gave some answers for some special case of the third question given by himself in \cite{Ho1}.
	
	 The second question given by Hong, which asked about the behaviour of the composition on the computation of reduced Gr\"obner bases,   has a complete answer given  in 1998 by J. Guti\'errez and R. Rubio San Miguel in   \cite{GR}. In Section~\ref{Reduced Gr\"obner Bases Under Composition}, we will study and comment on this  work. Besides this answer, they studied the same behaviour for the minimal Gr\"obner bases and gave a full answer for it. Also, they studied the special case studied by Hong in \cite{Ho2} for the reduced and minimal cases.
	
	 We will work and comment in Section~\ref{Gr\"obner Bases Under Composition: Different Polynomial Rings and Monomial Orderings} on the paper  \cite{LZW} published by Z. Liu and M. Wang in 2001. In this paper, they studied all the work described in all previous sections, but in a general case where we have different monomial orderings and polynomial rings before and after the composition.
	In Section~\ref{Homogeneous Gr\"obner Bases Under Composition}, the case of homogeneous Gr\"obner Bases will be studied using the paper  \cite{LW1}   published by J. Liu and M. Wang in 2006. This work is very important for our definition since it will be over homogeneous polynomials.
	Finally, in Section~\ref{Further results on Gr\"obner Bases Under Composition}, we will study the second paper \cite{LW2} of J. Liu and M. Wang in this subject published in 2007. They generalized in this paper the first work of Hong in \cite{Ho2} and the work done by themselves in \cite{LW1}. The studied case was the case of $\varGamma$-homogeneous Gr\"obner bases where special values of $\varGamma$ bring us back to the two previous works.

\section{Basic Definitions and Results}\label{S:Basic}
As we said before, this section contains some definitions and results of the theory of usual Gr\"obner bases. Note that most of the definitions and results given
 in this section was taken from the book of Cox,  Little and O'Shea \cite{CLO2}.

\begin{definition}
	A monomial ordering (term ordering) $<$ on the polynomial ring $K[x_1,\dots, x_n]$ is a relation $<$ on $\mathbb{Z}^n_{\geq0}$, or	equivalently, a relation on the set of monomials $\{x^\alpha : \alpha \in \mathbb{Z}^n_{\geq0}\}$, satisfying
	\begin{enumerate}
		\item $<$ is a total (or linear) ordering on $\mathbb{Z}^n_{\geq0}$.
		\item If $ \alpha < \beta$ and $\gamma \in \mathbb{Z}^n_{\geq0}$, then $\alpha +\gamma < \beta + \gamma$.                             		
		\item $<$ is a well-ordering on $\mathbb{Z}^n_{\geq0}$. This means that every non-empty subset of $\mathbb{Z}^n_{\geq0}$
		has a smallest element under $<$. In other words, if $A \subseteq \mathbb{Z}^n_{\geq0}$ is non-empty, then
		there is $\alpha \in A $ such that $\alpha < \beta$ for every $\beta \neq \alpha $ in $A$.
	\end{enumerate}
	
	Given a monomial ordering $<$, we say that $\alpha \leq \beta$ when either $\alpha < \beta$ or $\alpha = \beta$.
	
\end{definition}

\begin{definition}
	
	Let $f=\sum_{\alpha}a_{\alpha}x^{\alpha}$ be a non-zero polynomial in the polynomial ring $K[x_1, \dots, x_n]$ and let $<$
	be a monomial order.
	\begin{enumerate}
		\item  The multi-degree of $f$ with respect to $<$ is
		$$	multideg_{<}( f) = \max \{\alpha \in \mathbb{Z}^n_{\geq 0}: a_{\alpha}\neq 0 \}$$
		(the maximum is taken with respect to $<$).
		
		\item  The leading coefficient of $f$ with respect to $<$ is
		$$	LC_{<}( f) =  a_{multideg_{<}( f)} \in K.$$

		\item  The leading monomial of $f$ with respect to $<$ is
		$$ LM_{<}( f) = x^{multideg_{<}( f)} $$
		(with coefficient $1$).
		
		\item  The leading term of $f$ with respect to $<$ is
		$$	LT_{<}( f) = LC_{<}( f)\cdot LM_{<}( f).$$
	\end{enumerate}	
\end{definition}

\begin{lemma}
	Let $f, g \in K[x_1, \dots, x_n]$ be non-zero polynomials and let $<$ be any monomial ordering on $\mathbb{Z}^n_{\geq 0}$. Then
	\begin{enumerate}[leftmargin=15pt]
		\item $multideg_{<}( fg) = multideg_{<}( f) + multideg_{<}(g)$.
		\item  If $f + g \neq 0$, then $multideg_{<}( f + g) \leq \max (multideg_{<}( f ),multideg_{<}(g))$. If, in addition, $multideg_{<}( f ) \neq multideg_{<}(g)$, then equality occurs.
		\item $LT_{<}(f\cdot g) = LT_{<}(f) \cdot LT_{<}(g)$.
		\item $LM_{<}(f\cdot g) = LM_{<}(f) \cdot LM_{<}(g)$.
		\item If $LM_{<}(f) < LM_{<}(g)$, then $LT_{<}(f + g) = LT_{<}(g)$.
		\item If $LM_{<}(f) < LM_{<}(g)$, then $LM_{<}(f + g) = LM_{<}(g)$.
	\end{enumerate}
\end{lemma}

\begin{definition}
	
	Let $I \subseteq  K[x_1, \dots, x_n]$ be an ideal other than $\{0\}$, and fix a monomial
	ordering $<$ on $K[x_1, \dots, x_n]$.
	
	\begin{enumerate}
		\item  We denote by $LT_{<}(I)$ the set of leading terms of non-zero elements of $I$ with respect to $<$. Thus,
		$$LT_{<}(I) = \{cx^{\alpha} : \ \exists f \in I - \{0\} ~s.t.~ LT_{<}( f) = cx^{\alpha}\}.$$
		
		\item  We denote by $ \langle LT_{<}(I) \rangle$ the ideal generated by the elements of $LT_{<}(I)$.	
	\end{enumerate}	
\end{definition}

\begin{definition} \label{Gr\"obner basis}
	
	Fix a monomial ordering $<$ on the polynomial ring $K[x_1, \dots, x_n]$. A finite
	subset $G = \{g_1, \dots, g_t\}$ of non-zero polynomials of an ideal $I \subseteq  K[x_1, \dots, x_n]$ different from $\{0\}$ is said to
	be a Gr\"obner basis (or standard basis) with respect to $<$ if $$ \langle LT_{<}(I) \rangle = \langle LT_{<}(g_1), \dots, LT_{<}(g_t) \rangle.$$
	
	Using the convention that $ \langle \emptyset \rangle =  \{0\}$, we define the empty set $\emptyset $ to be the Gr\"obner
	basis for the zero ideal $\{0\}$. We will use the notation $G. B. _{<} (G, I)$ to say that $G$ is a Gr\"obner basis  with respect to $<$ for the ideal $I$.
\end{definition}

Note that if $F$ is not an ideal, then $G. B. _{<} (G, F)$ means that $G$ is a Gr\"obner basis  with respect to $<$ for the ideal $\langle F \rangle$. The same remark will be noted for any other type of Gr\"obner bases.

\begin{remark}
	Note that in the definition of a Gr\"obner basis $G$ with respect to some monomial ordering $ <$ for an ideal $I$, there are two different methods. The first which we used in Definition~\ref{Gr\"obner basis} and found in the book    \cite{CLO2} do not put that $G$ is a basis of $I$ as a condition in the definition. They put this fact as a result of the definition. The other way which found in the papers of Hong and other authors who studied the several types of Gr\"obner bases under composition (\cite{Ho1,Ho2,LW1,LW2,LZW}) put this result as a condition in their definitions.   	
\end{remark}

\begin{definition}
Let $f, g \in  K[x_1, \dots, x_n]$ be non-zero polynomials and fix a monomial
	ordering $<$ on $K[x_1, \dots, x_n]$.
	\begin{enumerate}[leftmargin=18pt]
		\item If $multideg_{<}( f) = \alpha $ and $multideg_{<}(g) = \beta$, then let $\gamma = (\gamma_1, \dots, \gamma_n)$, where $\gamma_i = \max(\alpha_i, \beta_i)$ for each i. We call $x^{\gamma}$ the least common multiple of $LM_{<}( f )$
		and $LM_{<}(g)$, written $x^{\gamma} = \lcm(LM_{<}( f ), LM_{<}(g))$. We will use later the following notation: $$ \sigma_<(f, g) = \dfrac{x^{\gamma}}{LT_{<}( f )} = \dfrac{\lcm(LM_{<}( f ), LM_{<}(g))}{LT_{<}( f )}.$$
		
		\item The $S_{<}$-polynomial of $f$ and $g$ with respect to $<$ is the combination
		$$S_{<}( f, g) = \dfrac{x^{\gamma}}{LT_{<}( f )}\cdot f-\dfrac{x^{\gamma}}{LT_{<}( g )}\cdot g = \sigma_<(f, g)\cdot f- \sigma_<(g, f)\cdot g.$$
	\end{enumerate}
\end{definition}

\begin{theorem}
	(Buchberger's Criterion). Let $I$ be a polynomial ideal and fix a monomial ordering $<$. Then a basis $G = \{g_1, \dots, g_t\}$ of $I$ is a Gr\"obner basis with respect to $<$ of $I$ if and only if for all pairs $i\neq j$, the remainder on division of $S_{<}(g_i, g_j)$ by $G$ (listed in some order) respect to $<$ is zero.
\end{theorem}
Buchberger's Criterion has other equivalent forms which we will set it here because we will need them later in the second chapter.
\begin{theorem}\label{Buchberger's Criterion1}
	Let $I$ be a polynomial ideal and fix a monomial ordering $<$. Then a basis $G = \{g_1, \dots, g_t\}$ of $I$ is a Gr\"obner basis with respect to $<$ of $I$ if and only if for all $g_i, g_j \in G, i \neq j$, there exist $h_1, \dots, h_t$ such that
	
	1.  $S_{<}(g_i, g_j) = h_1g_1+ \dots +h_tg_t$.
	
	2.  For every $k$, either $h_k = 0$ or \\
	\hspace*{4cm}  $LM_{<}(h_k) \cdot LM_{<}(g_k)  < \lcm (LM_{<}(g_i), LM_{<}(g_j))$.
\end{theorem}
This is one of the key theorems in Gr\"obner basis theory. Note that the statement of the theorem, in particular condition (2), is slightly different from the one usually found in the literature, Buchberger \cite{BB}; Cox, Little and O'Shea \cite{CLO2}; and Becker and Weispfenning \cite{BeWe}, in that $LT_{<}(S_{<}(g_i, g_j))$ is usually used in place of $LT_{<}(h_k)\cdot LT_{<}(g_k)$ and $\leq$ in place of $<$. However, the proofs for both are essentially the same \cite{Ho1}.
\begin{corollary}
	Let $I$ be a polynomial ideal and fix a monomial ordering $<$. Then a basis $G = \{g_1, \dots, g_t\}$ of $I$ is a Gr\"obner basis with respect to $<$ of $I$ if and only if for all $g_i, g_j \in G, i \neq j$, there exist $h_1, \dots, h_t$ such that
	
	1. $S_{<}(g_i, g_j) = h_1g_1+ \dots +h_tg_t$.
	
	2. For every $k$, either $h_k = 0$ or \\ \hspace*{4cm} $LT_{<}(h_k)\cdot LT_{<}(g_k) < \lcm (LT_{<}(g_i), LT_{<}(g_j))$.
	
	3. For every $k <l $, no term in $h_l\cdot LT_{<}(g_l)$ is divisible by $LT_{<}(g_k)$.
\end{corollary}
This is almost the same as in the previous theorem, except that we have one more condition (3). The first direction, in particular (3), follows immediately from the characterization of the generalized division described in \cite[p. 64]{CLO2}. The second direction is immediate from the previous theorem.	

\begin{lemma}
	Let $G$ be a Gr\"obner basis for $I \subseteq  K[x_1, \dots, x_n]$  with respect to $<$. Let $p \in G $ be a polynomial such that $LT_{<}(p) \in  \langle LT_{<}(G - \{p\})  \rangle$. Then $G - \{p\}$ is also a Gr\"obner basis for $I$ with respect to $<$.
\end{lemma}

\begin{definition}
	A minimal Gr\"obner basis for a polynomial ideal $I$ with respect to $<$ is a Gr\"obner basis
	$G$ for $I$ with respect to $<$ such that
	\begin{enumerate}
		\item $LC_{<}(p) = 1$ for all $p \in G$.
		
		\item For all $p \in G$, $LT_{<}(p) \notin\langle LT_{<}(G - \{p\})  \rangle$.
	\end{enumerate}	
	We will use the notation $M. G. B. _{<} (G, I)$ to say that $G$ is a minimal Gr\"obner basis  with respect to $<$ for the ideal $I$.
\end{definition}

\begin{definition}
	A reduced Gr\"obner basis for a polynomial ideal $I$ with respect to $<$ is a Gr\"obner basis
	$G$ for $I$ with respect to $<$ such that
	\begin{enumerate}
		\item $LC_{<}(p) = 1$ for all $p \in G$.
		
		\item For all $p \in G$, no monomial of $p$ lies in  $\langle LT_{<}(G - \{p\})  \rangle$.
	\end{enumerate}	
	We will use the notation $R. G. B. _{<} (G, I)$ to say that $G$ is a reduced Gr\"obner basis  with respect to $<$ for the ideal $I$. Note that we can define the same concepts for any non-empty subset of non-zero polynomials $F$ in $K[x_1, \dots, x_n]$ which need not to be a Gr\"obner basis and say that $F$ is a reduced set of polynomials with respect to $<$ if
	\begin{enumerate}
		\item $LC_{<}(f) = 1$ for all $f \in F;$ and
		\item For all $f \in F$, no monomial of $f$ lies in  $\langle LT_{<}(F - \{f\})  \rangle$.
	\end{enumerate}
	We denote a reduced set of polynomials $F$ with respect to $<$  by $R._{<} (F)$.	
\end{definition}

Note that for any polynomial ideal $I \neq \{0\} $and a given monomial ordering $<$, $I$ has a reduced Gr\"obner basis with respect to $<$, and this reduced Gr\"obner basis is unique.

Now, we will talk about the operation of composition of polynomials with several variables. We will give the definitions and set the previous results when we make the composition in the usual definition of leading terms and monomials. We will write and prove some similar results in the new definition also. This operation is very important and has a lot of applications. We take this paragraph from the paper of Hong  \cite{Ho1}:

``\emph{Composed objects (polynomials) often occur in real-life problem-solving because the underlying mathematical models are usually hierarchically structured. For instance, numerous physical quantities (such as work, torque, etc.) are defined in terms of other more basic quantities (such as length, time, etc.). Thus, we often need to deal with a set of polynomials in which the variables are defined in terms of other variables.}''

We need to use the following notations:
\begin{enumerate}
	\item $\Theta := (\theta_1, \dots, \theta_n)$, a list of $n$ non-zero polynomials in $K[x_1, \dots, x_n]$.
	\item $LT_{<}(\Theta) := (LT_{<}(\theta_1), \dots,LT_{<} (\theta_n))$.
	\item $LM_{<}(\Theta) := (LM_{<}(\theta_1), \dots,LM_{<} (\theta_n))$.
	\item $LC_{<}(\Theta) := (LC_{<}(\theta_1), \dots,LC_{<} (\theta_n))$.
	\item $\deg_{x_i}(f):=$ ~~the degree of $f$ with respect to the variable $x_i$.
\end{enumerate}
\begin{definition}[\cite{Ho1}] \label{def of Composition}
	(Composition) Let $h$ be a polynomial in the ring $K[x_1, \dots, x_n]$ and let $\Theta = (\theta_1, \dots, \theta_n)$ be a list of $n$ non-zero polynomials in $K[x_1, \dots, x_n]$. The composition of a polynomial $h$ by $\Theta$, written as $h \circ \Theta$, is the polynomial obtained from $h$ by replacing each $x_i$ in it with $\theta_i$. Likewise, $H \circ \Theta$ is the set $\{h \circ \Theta: h \in H\}$.
\end{definition}
One might consider the possibility of defining composition as the ``function composition'', namely, for all $(x_1, \dots, x_n ) \in K^n$, we have $$ (h\circ\Theta )(x_1,\dots, x_n ) = h( \theta_1(x_1,\dots, x_n ),\dots,\theta_n(x_1, \dots, x_n )).$$
But this is not suitable since $h \circ \Theta$ is not uniquely determined when $K$ is a finite field.

The composition by $\Theta$  is the endomorphism  $\Psi$ of $K[x_1, \dots, x_n]$
over $K$ that sends $x_i$ to $\theta_i$. The image of a polynomial $f$ w.r.t. such a composition is denoted by $f \circ \Theta$.
Similarly, for any set of polynomials $H \subset K[x_1, \dots, x_n]$, the set of images of elements of $H$ is denoted by
$H \circ \Theta$. If $I$ is an ideal of $K[x_1, \dots, x_n]$, the ideal generated by the image of $\Psi$, $ \langle \Psi(I) \rangle$, called extension ideal, is denoted by   $ \langle I \circ \Theta \rangle$.

Definition~\ref{def of Composition} can be generalized in the following way.
\begin{definition}[\cite{LZW}] \label{gene def of Composition}
	Assume that $n\leq m$ are two natural numbers. Let $h$ be a polynomial in the ring $K[x_1, \dots, x_n]$ and let $\Theta = (\theta_1, \dots, \theta_n)$ be a list of $n$ non-zero polynomials in the ring $K[y_1, \dots, y_m]$. The composition of a polynomial $h$ by $\Theta$, written as $h \circ \Theta$, is the polynomial obtained from $h$ by replacing each $x_i$ in it with $\theta_i$. Likewise, $H \circ \Theta$ is the set $\{h \circ \Theta: h \in H\}$.
\end{definition}
Note that in Definition~\ref{def of Composition}, $h \circ \Theta \in K[x_1, \dots, x_n]$, but in Definition~\ref{gene def of Composition} we have that $h \circ \Theta \in K[y_1, \dots, y_m]$. If $m=n$, then both definitions are the same. We will assume by default that $\Theta \in K[x_1, \dots, x_n]$ unless we say other things.

We set now two definitions which will be used later when we talk about the conditions on the list $\Theta$ which make the composition commute with the Gr\"obner basis computations.
\begin{definition}[\cite{Ho2}]
	Fix a monomial ordering $<$ on the polynomial ring $K[x_1, \dots, x_n]$. Let $\Theta$ be a list of $n$ non-zero polynomials in $K[x_1, \dots, x_n]$. The exponent matrix of $LM_{<}(\Theta)$ denoted by $\Mat(LM_{<}(\Theta))$, is the $n \times n$ matrix whose $(i, j)$-th entry is $\deg_{x_i}(LM_{<}(\theta_j))$. In other words, the $j$-th column of the matrix consists of the exponents of the leading monomial of $\theta_j$.
\end{definition}
\begin{definition}[\cite{Ho1}]
	(Permuted Powering) Fix a monomial ordering $<$ on the polynomial ring $K[x_1, \dots, x_n]$. Let $\Theta$ be a list of non-zero polynomials in $K[x_1, \dots, x_n]$. The list $LM_{<}(\Theta)$ is said to be permuted powering if  $$LM_{<}(\Theta)= (x^{\lambda_1}_{\pi_1},\dots,x^{\lambda_n}_{\pi_n})$$ for some permutation $\pi$ of $(1,\dots, n)$ and some $\lambda_1,\dots, \lambda_n > 0$.
\end{definition}

\begin{lemma} [\cite{Ho1}] \label{propert 1}
	Let $f, g \in K[x_1, \dots, x_n]$ be non-zero polynomials and $p \in K[x_1, \dots, x_n]$ be any non-zero monomial and $\Theta$ a list of $n$ non-zero polynomials in $K[x_1, \dots, x_n]$. Fix a monomial ordering $<$. Then
	\begin{enumerate}
		\item $(f\cdot g)\circ\Theta = f\circ\Theta \cdot g\circ\Theta$.
		\item $(f + g)\circ\Theta = f\circ\Theta + g\circ\Theta$.
		\item $LT_{<}(p\circ\Theta) = p\circ (LT_{<}(\Theta))$.
		\item $LM_{<}(p\circ\Theta) = p\circ (LM_{<}(\Theta))$.
	\end{enumerate} 	
\end{lemma}


\section{Gr\"obner Bases under Composition: Sufficiency}\label{Gr\"obner Bases Under Composition1}
The main question  is when the composition by some list of polynomials $\Theta $ brings a Gr\"obner basis with respect to some monomial ordering $<$ to a Gr\"obner basis with respect to the same ordering $<$ and the main contribution is to provide a simple answer to this question.
\begin{definition}[\cite{Ho1}] \label{Commutativity with Composition}
	(Commutativity with Composition) Fix a monomial ordering $<$ on the polynomial ring $K[x_1, \dots, x_n]$. We say that the composition by $\Theta $ commutes with the computation of Gr\"obner bases with respect to  $<$  if the following formula is true for $\Theta$:
 $$ \ \forall I \ \forall G \ [G. B._< (G, I) \Rightarrow G. B._< (G\circ\Theta, \langle I\circ\Theta \rangle)].$$
\end{definition}
It is easy to show that not every $\Theta$ satisfies this definition and we can easily construct counterexamples (for instance,
just permute the variables). But one can also find numerous positive examples. We also  want here to study the case when the list $\Theta$ maps every non-trivial Gr\"obner basis to a non-trivial one. This is very important and Hong used this thing during his study in the first paper \cite{Ho1} of them in the proof of the necessity part without talking about it directly in this definition and we will assume this fact in every section of our paper. As we said before in the introduction, Hong answered the question of this section and proved that the computation of Gr\"obner bases with respect to $<$ commutes with the composition if and only if the composition is ``compatible'' with the term ordering and the non-divisibility.
The definition of the compatibility with the term ordering and the non-divisibility will be given now.

\begin{definition}[\cite{Ho1}] \label{Compatibility with monomial ordering}
	(Compatibility with monomial ordering) Fix a monomial ordering $<$ on the polynomial ring $K[x_1, \dots, x_n]$. We say that the composition by $\Theta $ is compatible with the monomial ordering $<$  if  for all monomials $p$ and $q$, the following formula is true for $\Theta $:
 $$ \ \forall p \ \forall q \ [p < q \Rightarrow p\circ (LM_{<}(\Theta)) < q\circ (LM_{<}(\Theta))].$$
\end{definition}
\begin{definition}[\cite{Ho1}]
	(Compatibility of $<$ with non-divisibility) Fix a monomial ordering $<$ on the polynomial ring $K[x_1, \dots, x_n]$. We say that the composition by $\Theta $ is compatible with the non-divisibility with respect to $<$  if for all monomials $p$ and $q$, the following formula is true for  $\Theta $:
$$ \ \forall p \ \forall q \ [p \nmid q \Rightarrow p\circ (LM_{<}(\Theta)) \nmid q\circ (LM_{<}(\Theta))].$$
\end{definition}
The reason for using non-divisibility instead of divisibility is because  divisibility is compatible with every composition and thus not a useful condition. The following theorem is the main result of Sections~\ref{S:Basic} and  \ref{Gr\"obner Bases Under Composition1}.
\begin{theorem}[\cite{Ho1}] \label{Gr\"obner Bases Under Composition, main result}
	(Main theorem of commutativity with composition  with respect to $<$ ) Fix a monomial ordering $<$ on the polynomial ring $K[x_1, \dots, x_n]$. Let $\Theta$ be a list of non-zero polynomials in $K[x_1, \dots, x_n]$, then the following are equivalent
	\begin{itemize}
		\item (A) Composition by $\Theta $ commutes with the computation of Gr\"obner bases with respect to  $<$; and
		\item (B) Composition by $\Theta $ is
		\begin{enumerate}
			\item compatible with the monomial ordering $<$;  and
			\item compatible with the non-divisibility with respect to $<$.
		\end{enumerate}
	\end{itemize}	
\end{theorem}
In this section, it will be shown the sufficiency of this theorem, i.e. (B)  implies (A). To prove the two directions of Theorem~\ref{Gr\"obner Bases Under Composition, main result}, we must first set many lemmas which will describe the behaviour of composition when one or both of the two conditions of compatibility with monomial ordering or with the non-divisibility with respect to $<$ satisfied. We start by the following lemma, which states that a composition operation commutes with the leading monomial (term) extraction if it is compatible with the monomial ordering. In other words, how we can compute the leading term and monomial for the polynomial $f\circ\Theta$ for any polynomial $f$ if the composition by $\Theta $ is compatible with the monomial ordering $<$. The third part of Lemma~\ref{TO give 1}(B) will be used in Section~\ref{Reduced Gr\"obner Bases Under Composition}.
\begin{lemma} (\cite{GR,Ho1}) \label{TO give 1}
	Fix a monomial ordering $<$ on the polynomial ring $K[x_1, \dots, x_n]$, and let
	\begin{itemize}
		\item (A) Composition by $\Theta $ is compatible with the monomial ordering $<$;  and
		\item (B) For every $f \in K[x_1, \dots, x_n]$, we have
		\begin{enumerate}
			\item $LT_{<}(f\circ\Theta) = LT_{<}(f)\circ LT_{<}(\Theta)$.
			\item $LM_{<}(f\circ\Theta) = LM_{<}(f)\circ LM_{<}(\Theta)$.
			\item If $LC_{<}(\Theta)= (1, \dots,1)$ and $LC_{<}(f)=1$, then $LC_{<}(f\circ\Theta) = 1$.
		\end{enumerate}
	\end{itemize}
	Then (A)$\Rightarrow$ (B).	
\end{lemma}
The following lemma completely characterizes the condition of compatibility with the non-divisibility with respect to $<$ by giving an equivalent conditions of it.
These equivalent conditions will depend on the list $LM_{<}(\Theta)$ that defines the matrix $\Mat(LM_{<}(\Theta))$ which will be used in the proof.
The effects of this lemma will appear in the proof of the next lemmas and also in the construction of examples.
Sometimes, especially in the works  after the two papers of Hong (\!\!\cite{Ho1,Ho2}), whose studied other types of Gr\"obner bases under composition, the condition of permuted powering will be used instead of the compatibility with the non-divisibility. It will also be used in the next section while proving the necessity part.
\begin{lemma}[\cite{Ho1}] \label{ND give 1}
	Fix a monomial ordering $<$ on the polynomial ring $K[x_1, \dots, x_n]$, and let
	\begin{itemize}
		\item (A) Composition by $\Theta $ is compatible with the non-divisibility with respect to $<$;  and
		\item (B) the list $LM_{<}(\Theta)$ is permuted powering.
	\end{itemize}
	Then $(A)\Longleftrightarrow (B)$.	
\end{lemma}
 This lemma states that the composition  commutes with the least common multiple computation if it is compatible with the non-divisibility. Remember that the second condition of Theorem~\ref{Buchberger's Criterion1} depends on the least common multiple of the leading monomials of two polynomials and this theorem will by a very important key when we determine weather a given set is a Gr\"obner basis or not.
\begin{lemma}[\cite{Ho1}] \label{ND give 2}
	Fix a monomial ordering $<$ on the polynomial ring $K[x_1, \dots, x_n]$, and let $p$ and $q$ be any two monomials in $K[x_1, \dots, x_n]$. Let
	\begin{itemize}
		\item (A) Composition by $\Theta $ is compatible with the non-divisibility with respect to $<$;  and
		\item (B) $ \ \forall p \ \forall q [\lcm (p\circ (LM_{<}(\Theta)),q\circ (LM_{<}(\Theta)))  = \lcm (p,q) \circ (LM_{<}(\Theta))]$.
	\end{itemize}
	Then $(A)\Rightarrow (B)$.	
\end{lemma}
We will know from the next lemma that the operation of composition will preserve the generated ideals, i.e. if two different sets of polynomials generate the same ideal, then after composition they will still generate the same ideal.
\begin{lemma}[\cite{Ho1}]
	For any two subsets $G$ and $F$ of the polynomial ring $K[x_1, \dots, x_n]$ and any list $\Theta$ of non-zero polynomials in $K[x_1, \dots, x_n]$, we have
 $$ \langle G\rangle = \langle F \rangle\Rightarrow \langle G\circ\Theta \rangle =  \langle F\circ\Theta\rangle.$$
\end{lemma}
Lemma~\ref{comp gine not comp} will give us another form of the condition of the commutativity of the composition by $\Theta$. This form ignores the existence of the ideal $I$ in the original form and depends on the finite generating set $G$ directly. This is good for the calculations in the next results even in the next section also. Some authors used this form instead of the original one in their definitions, but we choose the original form because the first paper in our subject used it.
\begin{lemma}[\cite{Ho1}] \label{comp gine not comp}
	Let $G$ be any finite subset and $I$ be any ideal in the polynomial ring $K[x_1, \dots, x_n]$. Let $\Theta$ be a list of non-zero polynomials in $K[x_1, \dots, x_n]$. Fix a monomial ordering $<$ on $K[x_1, \dots, x_n]$. Let
	\begin{itemize}
		\item (A) $ \ \forall I \ \forall G \ [G. B._< (G, I) \Rightarrow G. B._< (G\circ\Theta, \langle I\circ\Theta \rangle)];$ and
		\item (B) $ \ \forall G \ [G. B._< (G) \Rightarrow G. B._< (G\circ\Theta)]$.
	\end{itemize}
	Then $(A)\Longleftrightarrow (B)$.
\end{lemma}
The next lemma is the key of the sufficiency part of the main theorem of these two sections and all previous results were written for it. The second form of the condition of the commutativity of the composition by $\Theta$ will be used instead of the original one and in Theorem~\ref{Gr\"obner Bases Under Composition, suff} we will back to it. In the proof, Theorem~\ref{Buchberger's Criterion1} will be used to test weather a given set (before and after the composition) is a Gr\"obner basis with respect to some monomial ordering $<$ or not. The effects of compatibility with the monomial ordering $<$ on the computation of the leading terms and monomials after composition and the effects of compatibility with the non-divisibility with respect to $<$ (or the permuted powering condition) on the computing least common multiple after the composition will be used also.
\begin{lemma}[\cite{Ho1}] \label{not comp, suff}
	Fix a monomial ordering $<$ on the polynomial ring $K[x_1, \dots, x_n]$, and let $G$ be any subset of $K[x_1, \dots, x_n]$. Let $\Theta$ be a list of non-zero polynomials in $K[x_1, \dots, x_n]$. Let
	\begin{itemize}
		\item (A) $ \ \forall G [G. B._< (G) \Rightarrow G. B._< (G\circ\Theta)];$ and
		\item (B) Composition by $\Theta $ is
		\begin{enumerate}
			\item compatible with the monomial ordering $<$; and
			\item compatible with the non-divisibility with respect to $<$.
		\end{enumerate}	
	\end{itemize}
	Then $(B)\Rightarrow (A)$.	
\end{lemma}
Finally, we are ready to state the sufficiency of the main theorem. The difference between this theorem and Lemma~\ref{not comp, suff} is only in the forms of the condition of the commutativity of the composition by $\Theta$.
\begin{theorem}[\cite{Ho1}] \label{Gr\"obner Bases Under Composition, suff}
	(Main theorem of commutativity with composition  with respect to $<$, Sufficiency) Fix a monomial ordering $<$ on the polynomial ring $K[x_1, \dots, x_n]$. Let $\Theta$ be a list of non-zero polynomials in $K[x_1, \dots, x_n]$. Let
	\begin{itemize}
		\item (A) Composition by $\Theta $ commutes with the computation of Gr\"obner bases with respect to  $<$; and
		\item (B) Composition by $\Theta $ is
		\begin{enumerate}
			\item compatible with the monomial ordering $<$; and
			\item compatible with the non-divisibility with respect to $<$.
		\end{enumerate}
	\end{itemize}
	Then $(B)\Rightarrow (A)$.		
\end{theorem}
\section{Gr\"obner Bases under Composition: Necessity}\label{Gr\"obner Bases Under Composition2}
In this section, the proof of the necessity part of the main theorem of the compatibility conditions for commutativity will be described, i.e.  (A)  implies (B). There is no  way to describe the strategy of the proof better than what Hong, who invented  it, said in his paper \cite{Ho1}:

``\emph{Before plunging into the detail of the `long' proof, we describe the overall strategy. Mostly the proof is by proving contrapositive. Thus, it goes like this. Assume that (B) is not true. Then find $G$ such that $G. B._<(G)$ but not $G. B._<(G\circ\Theta)$. Obviously the main difficulty in this process lies in finding such $G$. I had to spend numerous days (experimenting with
	computer algebra systems, making conjectures, disproving them to my dismay, dreaming about them in my sleep, etc., as usual) to find the ones presented here. Once they have been found, it was easy to write down the `straight-line forward' proof. Lemmas \ref{not comp give TO} and \ref{not gb1} are the cores of the proof, that is, they contain such $G$s as those mentioned above.}''

After Hong, every author studied the behaviour of composition over several types of Gr\"obner bases in the usual sense including our study over fields with valuations used this strategy in their proofs of the necessity parts. Note that, as Hong did in his paper and we explained in the previous section, the second form of the condition of the commutativity of the composition by $\Theta$ will be used instead of the original one in all results in this section.

Now, to show that commutativity of the composition implies the compatibility with the monomial ordering, we need the following three lemmas where:
\begin{itemize}
	\item In Lemma~\ref{not comp give  noneq1}, non-equal monomials composite to non-equal polynomials and this non-equality does not depend on a constant product. The first non-equality comes from
$$~p < q ~\text{and}~a\neq 0 ~\text{and}~b \neq 0~\Rightarrow ap \neq bq.$$
 Note also we need here the fact that we work with the list $\Theta$ which maps every non-trivial ideal to a non-trivial one. If we do not set this condition, then this lemma will be false.
	\item In Lemma~\ref{not comp give  noneq2}, non-equal monomials composite to two polynomials with non-equal leading monomials. The proof of it depends on the result proved in Lemma~\ref{not comp give  noneq1} and the fact that $LM_{<}( p\circ \Theta) = p\circ LM_{<}(\Theta)$.
	\item In Lemma~\ref{not comp give TO}, we finish to needed result, that is the composition by $\Theta$ compatible with the monomial ordering $<$. As above, the proof of it depends on Lemma~\ref{not comp give  noneq2}.
\end{itemize}
\begin{lemma}[\cite{Ho1}] \label{not comp give  noneq1}
	Let $G$ be any subset of  $K[x_1, \dots, x_n]$ and let $p$ and $q$ be any two monomials in $K[x_1, \dots, x_n]$ and let $a, b \in K$. Let $\Theta$ be a list of non-zero polynomials in $K[x_1, \dots, x_n]$. Fix a monomial ordering $<$ on $K[x_1, \dots, x_n]$. Let
	\begin{itemize}
		\item (A) $  \ \forall G  \ [G. B._< (G) \Rightarrow G. B._< (G\circ\Theta)];$ and
		\item (B) $ \ \forall p \ \forall q \ \forall a \ \forall b  \ [~ p < q ~\text{and}~a\neq 0 ~\text{and}~b \neq 0 ~ \Rightarrow  ap\circ\Theta \neq bq\circ\Theta]$.
	\end{itemize}
	Then (A) $\Rightarrow$ (B).
\end{lemma}
\begin{lemma}[\cite{Ho1}] \label{not comp give  noneq2}
	Let $G$ be any subset of  $K[x_1, \dots, x_n]$ and let $p$ and $q$ be any two monomials in $K[x_1, \dots, x_n]$. Let $\Theta$ be a list of non-zero polynomials in $K[x_1, \dots, x_n]$. Fix a monomial ordering $<$ on $K[x_1, \dots, x_n]$. Let
	\begin{itemize}
		\item (A) $ \ \forall G \ [G. B._< (G) \Rightarrow G. B._< (G\circ\Theta)];$ and
		\item (B) $ \ \forall p \ \forall q \ [~ p < q  \Rightarrow p\circ (LM_{<}(\Theta))\neq q\circ (LM_{<}(\Theta))]$.
	\end{itemize}
	Then (A) $\Rightarrow$ (B).
\end{lemma}
\begin{lemma}[\cite{Ho1}] \label{not comp give TO}
	Let $G$ be any subset of $K[x_1, \dots, x_n]$ and let $p$ and $q$ be any two monomials in $K[x_1, \dots, x_n]$. Let $\Theta$ be a list of non-zero polynomials in $K[x_1, \dots, x_n]$. Fix a monomial ordering $<$ on $K[x_1, \dots, x_n]$. Let
	\begin{itemize}
		\item (A) $ \ \forall G \ [G. B._< (G) \Rightarrow G. B._< (G\circ\Theta)]$;
		\item (B) $ \ \forall p \ \forall q  \ [~ p < q  \Rightarrow p\circ (LM_{<}(\Theta)) < q\circ (LM_{<}(\Theta))]$.
	\end{itemize}
	Then (A) $\Rightarrow$ (B).
\end{lemma}
Thus, the proof of one half is finished; the commutativity implies the compatibility with the monomial ordering. Now, let us work on the other half: the commutativity implies the compatibility with the non-divisibility. We start by the following two lemmas where the second lemma looks like a generalization of the first one and plays an important role and will be the core of proving the compatibility with the non-divisibility.
They also explain how the non-relatively prime of the leading monomials with respect to some $<$ of two polynomials $f$ and $g$ affect on the existence of a Gr\"obner basis with respect to the same $<$ consisting of these two polynomials or some power of one and the same power of the other one plus $1$. Note that the non-relatively prime condition in the first one is not written directly.
\begin{lemma}[\!{\cite[Lemma 5.5]{Ho1}}] \label{not gb1}
	Fix a monomial ordering $<$ on the polynomial ring $K[x_1, \dots, x_n]$. Let $f$ and $g$ be  two non-zero polynomials in $K[x_1, \dots, x_n]$ with $LM_{<}(f)= x_1^{u_1}\dots x_n^{u_n}$ and $LM_{<}(g)= x_1^{v_1}\dots x_n^{v_n}$. Let $\Theta$ be a list of non-zero polynomials in $K[x_1, \dots, x_n]$. Assume that $0 < v_k \leq u_k $ for some $k$, then we have
	\begin{enumerate}
		\item  $ \{f,g\}$ is not a Gr\"obner basis with respect to $<$, or
		\item  $ \{f+1,g\}$ is not a Gr\"obner basis with respect to $<$.
	\end{enumerate}
\end{lemma}

\begin{lemma}[\!{\cite[Lemma 5.6]{Ho1}}] \label{not gb2}
	Fix a monomial ordering $<$ on the polynomial ring  $K[x_1, \dots, x_n]$. Let $f$ and $g$ be  two non-zero polynomials in $K[x_1, \dots, x_n]$ with $LM_{<}(f)= x_1^{u_1}\dots x_n^{u_n}$ and $LM_{<}(g)= x_1^{v_1}\dots x_n^{v_n}$. Let $\Theta$ be a list of non-zero polynomials in $K[x_1, \dots, x_n]$. Assume that the leading terms are not relatively prime, that is $0 < v_k $ and $0 < u_k $ for some $k$. Then there exists $\lambda>0$ such that
	\begin{enumerate}
		\item  $ \{f^\lambda,g\}$ is not a Gr\"obner basis with respect to $<$, or
		\item  $ \{f^\lambda+1,g\}$ is not a Gr\"obner basis with respect to $<$.
	\end{enumerate}
\end{lemma}
Now we work with the list of leading monomials $LM_{<}(\Theta)$ of the list $\Theta$ and describe how we can prove that this list is permuted powering if the composition by $\Theta$ commutes with the computation of Gr\"obner bases with respect to $<$. Assuming that this is true, then:
\begin{itemize}
	\item In Lemma~\ref{not comp give rp}, the list $LM_{<}(\Theta)$ consists of a pair-wise relatively prime monomials. We need here the result proved in Lemma~\ref{not gb2} for two polynomials with a non-relatively prime leading monomial.
	
	\item In Lemma~\ref{not comp give not const}, the result that the commutativity implies the compatibility with the monomial ordering given in Lemma~\ref{not comp give TO} is used to show that the list $LM_{<}(\Theta)$ does not contain any non constant monomials and that $LM_{<}(\theta_i)\neq 1$ for any $i$.
	
	\item In the proof of  Lemma~\ref{not comp give pp}  both Lemmas~\ref{not comp give rp} and \ref{not comp give not const} will be used.
	From Lemma~\ref{not comp give TO}, we know that the list $LM_{<}(\Theta)$ consists of a pair-wise relatively prime monomials.
	Therefore there exists at most one non-zero element in each row of the matrix $\Mat(LM_{<}(\Theta))$. From Lemma~\ref{not comp give not const},
	we also know that $LM_{<}(\theta_i)\neq 1$ for any $i$. Therefore there exists at least one non-zero element in each column of $\Mat(LM_{<}(\Theta))$. Thus, we see that there is exactly one non-zero element in each row and each column
	of $\Mat(LM_{<}(\Theta))$. Hence $\Mat(LM_{<}(\Theta))$ is a permuted diagonal matrix, which is equivalent to that list $LM_{<}(\Theta)$ is permuted powering.
	
	\item In Lemma~\ref{not comp give nd}, we finish to the needed result, that is the composition by $\Theta$ is compatible with the non-divisibility by replacing it by its equivalent condition, $LM_{<}(\Theta)$ is permuted powering, proved in Lemma~\ref{not comp give pp}.
\end{itemize}
\begin{lemma}[\cite{Ho1}] \label{not comp give rp}
	Let $G$ be any subset of  polynomials of the polynomial ring $K[x_1, \dots, x_n]$. Let $\Theta$ be a list of non-zero polynomials in $K[x_1, \dots, x_n]$. Fix a monomial ordering $<$ on $K[x_1, \dots, x_n]$. Let
	\begin{itemize}
		\item (A) $ \ \forall G \ [G. B._< (G) \Rightarrow G. B._< (G\circ\Theta)]$;
		\item (B) the monomials $LM_{<}(\theta_1), \dots,LM_{<} (\theta_n)$ be pair-wise relatively prime.
	\end{itemize}
	Then (A) $\Rightarrow$ (B).
\end{lemma}
\begin{lemma}[\cite{Ho1}] \label{not comp give not const}
	Let $G$ be any subset of polynomials of the polynomial ring $K[x_1, \dots, x_n]$. Let $\Theta$ be a list of non-zero polynomials in $K[x_1, \dots, x_n]$. Fix a monomial ordering $<$ on $K[x_1, \dots, x_n]$. Let
	\begin{itemize}
		\item (A) $ \ \forall G \ [G. B._< (G) \Rightarrow G. B._< (G\circ\Theta)]$;
		\item (B) $\ \forall i \ LM_{<}(\theta_i)\neq 1$.
	\end{itemize}
	Then (A) $\Rightarrow$ (B).
\end{lemma}
\begin{lemma}[\cite{Ho1}] \label{not comp give pp}
	Let $G$ be any subset of polynomials of the polynomial ring $K[x_1, \dots, x_n]$. Let $\Theta$ be a list of non-zero polynomials in $K[x_1, \dots, x_n]$. Fix a monomial ordering $<$ on $K[x_1, \dots, x_n]$. Let
	\begin{itemize}
		\item (A) $ \ \forall G \  [G. B._< (G) \Rightarrow G. B._< (G\circ\Theta)]$;
		\item (B) the list $LM_{<}(\Theta)$ is permuted powering.
	\end{itemize}
	Then (A) $\Rightarrow$ (B).
\end{lemma}
\begin{lemma}[\cite{Ho1}] \label{not comp give nd}
	Let $G$ be any subset of  polynomials of the polynomial ring $K[x_1, \dots, x_n]$. Let $\Theta$ be a list of non-zero polynomials in $K[x_1, \dots, x_n]$. Fix a monomial ordering $<$ on $K[x_1, \dots, x_n]$. Let
	\begin{itemize}
		\item (A) $ \ \forall G \ [G. B._< (G) \Rightarrow G. B._< (G\circ\Theta)]$;
		\item (B) $ \ \forall p \ \forall q \ [~p \nmid q \Rightarrow p\circ (LM_{<}(\Theta)) \nmid q\circ (LM_{<}(\Theta))]$.
	\end{itemize}
	Then (A) $\Rightarrow$ (B).
\end{lemma}
Finally, by joining the results proved in Lemma~\ref{not comp give TO} and Lemma~\ref{not comp give nd}, we get the next lemma which improved into Theorem \ref{Gr\"obner Bases Under Composition, nece} by going back to the original condition of the commutativity to get the necessity part of the main theorem.
\begin{lemma}[\cite{Ho1}]
	Fix a monomial ordering $<$ on the polynomial ring $K[x_1, \dots, x_n]$, and let $G$ be any subset of $K[x_1, \dots, x_n]$. Let $\Theta$ be a list of non-zero polynomials in $K[x_1, \dots, x_n]$. Let
	\begin{itemize}
		\item (A) $ \ \forall G \ [G. B._< (G) \Rightarrow G. B._< (G\circ\Theta)]$;
		\item (B) Composition by $\Theta $ is
		\begin{enumerate}
			\item compatible with the monomial ordering $<$; and
			\item compatible with the non-divisibility with respect to $<$.
		\end{enumerate}
	\end{itemize}	
	Then (A) $\Rightarrow$ (B).	
\end{lemma}
The last form of the necessity of the main theorem will be given now.
\begin{theorem}[\cite{Ho1}] \label{Gr\"obner Bases Under Composition, nece}
	(Main theorem of commutativity with composition  with respect to $<$, Necessity) Fix a monomial ordering $<$ on the polynomial ring $K[x_1, \dots, x_n]$. Let $\Theta$ be a list of non-zero polynomials in $K[x_1, \dots, x_n]$. Let
	\begin{itemize}
		\item (A) Composition by $\Theta $ commutes with the computation of Gr\"obner bases with respect to  $<$; and
		\item (B) Composition by $\Theta $ is
		\begin{enumerate}
			\item compatible with the monomial ordering $<$; and
			\item compatible with the non-divisibility with respect to $<$.
		\end{enumerate}
	\end{itemize}
	Then (A) $\Rightarrow$ (B).	
\end{theorem}
At the end of this study, we put some examples given by Hong at the end of his paper, starting from some trivial cases and finishing with a special monomial ordering.
\begin{example}[\cite{Ho1}] \label{ex1,1}
	Every composition of the form $$LM_<(\theta_i)=x_i^\lambda,$$ where $0 < \lambda $, is a compatible composition since it  satisfies the two compatibility conditions, and thus commutes with the computation of Gr\"obner bases with respect to $<$. Note that $<$ is arbitrary here but  $0 < \lambda $ is fixed and the permutation is the identity one. This mentioned class of composition covers many naturally arising
	compositions like:
	\begin{enumerate}
		\item Scaling: $\theta_i=a_ix_i,~a_i \neq 0$.\\ For example, $\Theta=(2x_1,3x_2)$.
		\item Translation: $\theta_i=x_i-c_i$.\\ For example, $\Theta=(x_1-2,x_2+3)$.
		\item Powering: $\theta_i=x_i^\lambda,~0 < \lambda$.\\ For example, $\Theta=(x_1^2,x_2^2)$.
		\item Univariate: $\theta_i \in K[x_i] $ of degree $0 < \lambda$.\\ For example, $\Theta=(2x_1^4 - 5x_1^2 +3x_1^2+ 4,5x_2^4 +2x_2^2 -3x_2-3)$.
		\item General: $\theta_i \in K[x_1, \dots, x_n]$ such that $LM_<(\theta_i)=x_i^\lambda,~0 < \lambda$. \\ For example, $\Theta=(2x_1^4 - 5x_1x_2^2 +4x_2^3+ 1,x_2^4 -2x_1^2x_2^2 +3x_12x_2-3)$, for the graded lexicographic ordering with $(x_1 < x_2)$.
	\end{enumerate}	
\end{example}
\begin{example}[\cite{Ho1}] \label{ex2,1}
	Let $<$ be the lexicographic ordering. Then, every composition of the form $$LM_<(\theta_i)=x_i^{\lambda_i},$$ where $0 < \lambda_i $, is a compatible composition
	and commutes with the computation of Gr\"obner bases with respect to $<$. Note that, different  $0 < \lambda_i $ allowed now for different  $x_i$. Some several compatible compositions for the lexicographic ordering will be listed below.
	\begin{enumerate}
		\item Powering: $\theta_i=x_i^{\lambda_i},~0 < \lambda_i$.\\ For example, $\Theta=(x_1^2,x_2^3)$.
		\item Univariate: $\theta_i \in K[x_i] $ of degree $0 < \lambda_i$.\\ For example, $\Theta=(2x_1^5 - 5x_1^2 +3x_1^2+ 4,5x_2^3 +2x_2^2 -3x_2-3)$.
		\item General: $\theta_i \in K[x_1, \dots, x_n]$ such that $LM_<(\theta_i)=x_i^{\lambda_i},~0 < \lambda_i$.\\ For example, $\Theta=(2x_1^3 - 5x_1+ 1,x_2^2 -2x_1^2x_2 +3x_1^5-3)$,  for the lexicographic ordering with $(x_1 < x_2)$.
	\end{enumerate}	
\end{example}
All the previous examples have one thing in common: $LM_<(\theta_i)$ involves $x_i$, that is, no permutation of variables. The following example is with a non-trivial permutation.
\begin{example}[\cite{Ho1}] \label{ex3,1}
	Let $p=x^{a_1}y^{a_2}$ and $q=x^{b_1}y^{b_2}$ be two monomials in $K[x,y]$. Consider the monomial ordering defined by
	$$ p < q \Longleftrightarrow a_1 +\sqrt{2} a_2 < b_1 +\sqrt{2} b_2.$$
Let $\Theta=(y + x, x^2 +y)$, then after calculations we have that $$LM_<(\Theta)=(y, x^2 ),$$ and thus the variables permute and $LM_<(\Theta)$ is permuted powering. For the compatibility with the monomial ordering, we have that if $ p < q$, then
	$$ p \circ LM_<(\Theta) = y^{a_1}x^{2a_2}= x^{2a_2} y^{a_1},~~q \circ LM_<(\Theta) = y^{b_1}x^{2b_2}= x^{2b_2} y^{b_1}.$$
 Thus,
	$$ 2a_2 +\sqrt{2} a_1 =\sqrt{2} (a_1 +\sqrt{2} a_2) < \sqrt{2} (b_1 +\sqrt{2} b_2) = 2b_2 +\sqrt{2} b_1,$$ hence $$ p \circ LM_<(\Theta) < q \circ LM_<(\Theta).$$
 This means that the composition by $\Theta=(y + x, x^2 +y)$ is a compatible composition and commutes with the computation of Gr\"obner bases with respect to $<$ defined above.
\end{example}
\section{Gr\"obner Bases under Composition: Different Monomial Orderings}\label{Gr\"obner Bases Under Composition3}
\sectionmark{Gr\"obner Bases under Composition: Different  Orderings}
In the last two previous sections we saw the answer of the question: when does the computation of Gr\"obner bases with respect to some monomial ordering $<$ commutes with composition? As we said before, the first one who studied this subject was Hong in  \cite{Ho1}. At the end of this paper, Hong said that this is not the end and gave three open questions to be answered. The third question was:

\emph{Let $G$ be a Gr\"obner basis for $I$ with respect to $<$. When is $G \circ \Theta$ is a Gr\"obner basis for $\langle I \circ \Theta \rangle$ (possibly with respect to another monomial  ordering $<'$)?}

He commented on this question by:

\emph{In order to answer this question, one could carefully analyse the proof given in this paper \cite{Ho1}, and generalize it. In fact, the author has already followed this approach and found some answer, which is reported in another paper \cite{Ho2}, but it might be interesting to find a completely new approach.}

There are two different meanings for the expression `another monomial ordering'.
\begin{itemize}
	\item The first meaning is to talk a bout `another \textit{fixed} monomial ordering'. That is, another monomial ordering may be different on the first one but has a special way to define it which may be using the first defined monomial ordering itself or using the list $\Theta$. In this section, we will study this type of understanding with a special way of defining the other ordering using both of the first fixed one and the list $\Theta$. We will see also what conditions  we must  have to make this new ordering be a monomial ordering. Note that some times the new obtained ordering will be equal to the first fixed one. In a mathematical writing, the problem which will be discussed in this section is which conditions must $\Theta$ satisfied to make the following statement true:
	\[ \ \forall I\ \forall G  \ [G. B._{\Theta\circ<}(G,I) \Rightarrow G. B._< (G\circ\Theta, \langle I\circ\Theta \rangle)],\] where $\Theta\circ<$ is some ordering defined later. The reader may be confused when he see that the obtained ordering will be used before the composition and the fixed one will be used after this operation. The reason for this order of using is because the list of leading monomials of the list $\Theta$ will be computed using the fixed ordering and the definition of the obtained one depends on this list. If we use the first ordering to define the Gr\"obner bases with respect to it before the composition, then we will have nothing to study and will back to the original question. This way of definition and the results for it were studied by Hong in  \cite{Ho2}.
	He proved that if the list of leading monomials with respect to some $<$ of the list $\Theta$ is permuted powering, then for any Gr\"obner basis $G$ of an ideal $I$ with respect to the obtained ordering (which will be monomial ordering under these conditions), the set  $G \circ \Theta$ is a Gr\"obner basis for $\langle I \circ \Theta \rangle$ with respect to $<$.
	\item The second meaning, which is the general one, is when both of the two monomial orderings are arbitrary without any special conditions on any of them. In other words, we choose two monomial orderings from the beginning and use one to define the Gr\"obner basis $G$ with respect to it before the composition and use the other ordering for $G \circ \Theta$ after this operation. This question and some other related question was studied by  Z. Liu and M. Wang in   \cite{LZW},  which we will study it in Section~\ref{Gr\"obner Bases Under Composition: Different Polynomial Rings and Monomial Orderings}.
\end{itemize}
We will start by the definition of the ordering $\Theta\circ<$ obtained by compose some monomial ordering $<$ with a list $\Theta$. Note that this ordering is a relation defined on the set of all monomials on the variables $\{x_1, \dots, x_n\}$. In general, not any pair of these monomials are related using $\Theta\circ<$. More explaining will be given later.
\begin{definition}[\cite{Ho2}] \label{Composition on Ordering}
	(Composition on Ordering) Fix a monomial ordering $<$ on the polynomial ring $K[x_1, \dots, x_n]$. Let $\Theta$ be a list of non-zero polynomials in $K[x_1, \dots, x_n]$. The composition of $<$ by $\Theta$, written as $\Theta\circ<$, is the binary relation over
	the monomials defined by $$\ \forall p \ \forall q \  [~p ~\Theta\circ< q \Longleftrightarrow p\circ (LM_{<}(\Theta)) < q\circ (LM_{<}(\Theta))].$$
\end{definition}
To make writing more easy, we will use the notation $<_\Theta$ instead of $\Theta\circ<$ in this paper.
Note that the relation $<_\Theta$ is not necessarily a monomial ordering. As a simple counterexample,
consider a univariate case where $\Theta = (1)$. Obviously $$1 \circ (LM_{<}(\Theta))= 1= x_1\circ (LM_{<}(\Theta)),$$ and thus it not true that $1 < x_1 $,
which violates one of the conditions of a monomial ordering. As we said before, the conditions which make it monomial ordering will be shown in Lemma~\ref{MO4}.
That is, the relation $<_\Theta$ is indeed a monomial ordering if and only if the exponent matrix $\Mat (LM_{<}(\Theta))$ is non-singular. Remember that from Lemma~\ref{ND give 1}, where we comment to the reader to look on its proof in Hong \cite{Ho1},
we have that the list $LM_{<}(\Theta)$ is permuted powering if and only if the exponent matrix $\Mat (LM_{<}(\Theta))$ is a permutation of a diagonal matrix with all positive diagonal entries if and only if it is non-singular.

Now to reach Lemma~\ref{MO4}, we must use the following three lemmas:
\begin{itemize}
	\item From Lemma~\ref{MO1}, we will know that the condition $$\ \forall p \ \forall q \ \forall r[~p <_\Theta q \Rightarrow ~rp <_\Theta rq]$$ is true for arbitrary $\Theta$. This happened since the monomial ordering $<$ which defines $<_\Theta$ satisfies this condition. Also, because the composition of polynomials with any $\Theta$ distributes under the polynomials product and $$LT_{<}(p\circ\Theta) = p\circ LT_{<}(\Theta)$$ for any monomial $p$ as we saw in Lemma~\ref{propert 1} before.
	
	\item In Lemma~\ref{MO2}, it will be shown that the second condition for monomial ordering (every non-trivial monomial must be greater than the trivial one) is equivalent to that every column of the exponent matrix $\Mat (LM_{<}(\Theta))$ has at least one non-zero entry. The proof  uses the concepts of linear algebra by showing that these two conditions are equivalent to the condition that the product of $\Mat (LM_{<}(\Theta))$ by any non-zero vector of $n$ natural numbers must be non-zero. Remember that if $\Mat (LM_{<}(\Theta))$ in non-singular, then the last statement must be true and this means that the second condition for monomial ordering will be satisfied.
	
	\item In Lemma~\ref{MO3}, the non-singularity of the matrix $\Mat (LM_{<}(\Theta))$ will be sufficient and necessary for the ordering $<_\Theta$ to be linear.
	Note $<_\Theta$ is clearly reflexive and the transitivity of it follows immediately from the transitivity of $<$ and does not need any other condition. Also, since for any two monomials $p\neq q$, we have that $p\circ LT_{<}(\Theta) < q\circ LT_{<}(\Theta)$ or $q\circ LT_{<}(\Theta) < p\circ LT_{<}(\Theta)$, then at most one of $p <_\Theta q$ or $q <_\Theta p$ must holds for any arbitrary $\Theta$.
	The non-singularity condition is used to prove that every two non-equal
	monomials $p\neq q$, must be comparable using $<_\Theta$ and equivalent to this fact. The core of the proof of this equivalent depends on the contrapositive of the fact that $\Mat (LM_{<}(\Theta))$ is non-singular if and only if $\det(\Mat (LM_{<}(\Theta))) \neq 0$.	
\end{itemize}

\begin{lemma}[\cite{Ho2}] \label{MO1}
	Fix a monomial ordering $<$ on the polynomial ring $K[x_1, \dots, x_n]$. Let $\Theta$ be a list of non-zero polynomials in $K[x_1, \dots, x_n]$. If $p, q$ and $r$ are monomials in $K[x_1, \dots, x_n]$, then we have $$\ \forall p \ \forall q \ \forall r \ [~p <_\Theta q \Rightarrow ~rp <_\Theta rq]$$is true for any arbitrary $\Theta$.
\end{lemma}

\begin{lemma}[\cite{Ho2}] \label{MO2}
	Fix a monomial ordering $<$ on the polynomial ring $K[x_1, \dots, x_n]$. Let $\Theta$ be a list of non-zero polynomials in $K[x_1, \dots, x_n]$ and $p$ is any monomial in $K[x_1, \dots, x_n]$. Let
	\begin{itemize}
		\item (A) Every column of the exponent matrix $\Mat (LM_{<}(\Theta))$ has at least one non-zero entry; and
		\item (B) $\ \forall p \ [~p \neq 1 \Rightarrow 1 <_\Theta p]$.
	\end{itemize}
	Then (A) $\Longleftrightarrow$ (B).
\end{lemma}

\begin{lemma}[\cite{Ho2}] \label{MO3}
	Fix a monomial ordering $<$ on the polynomial ring $K[x_1, \dots, x_n]$. Let $\Theta$ be a list of non-zero polynomials in $K[x_1, \dots, x_n]$. Let
	\begin{itemize}
		\item (A) The exponent matrix $\Mat (LM_{<}(\Theta))$ is non-singular; and
		\item (B) The binary relation $<_\Theta$ is a linear ordering.
	\end{itemize}
	Then (A) $\Longleftrightarrow$ (B).
\end{lemma}
Now, the last three lemmas will be joined together to form Lemma~\ref{MO4} which gives us the sufficient and necessary condition for the ordering $<_\Theta$ to be a monomial ordering.
\begin{lemma}[\cite{Ho2}] \label{MO4}
	Fix a monomial ordering $<$ on the polynomial ring $K[x_1, \dots, x_n]$. Let $\Theta$ be a list of non-zero polynomials in $K[x_1, \dots, x_n]$. Let
	\begin{itemize}
		\item (A) The exponent matrix $\Mat (LM_{<}(\Theta))$ is non-singular; and
		\item (B) The binary relation $<_\Theta$ is a monomial ordering.
	\end{itemize}
	Then (A) $\Longleftrightarrow$ (B).
\end{lemma}
From now, the main strategy of the proof is similar to what used in \cite{Ho2} by Hong itself, but with some different conditions. The reader will notice that the proof of the main result here is like a generalization of what given in \cite{Ho2}. Also, the reader should be careful with the essential differences, especially,  to the places where $<_\Theta$ is used instead of $<$. The core of the proof of the main theorem will be the results given in Lemma~\ref{MO7}, which is in fact the summary of the three Lemmas~\ref{MO4}, \ref{MO5} and \ref{MO6}, where:

\begin{itemize}
	\item In Lemma~\ref{MO5}, if the ordering $<_\Theta$ is a monomial ordering, then we can compute the leading term (monomial) with respect to $<$ for the polynomial $f\circ\Theta$ for any polynomial $f$ to be the composition  of the leading term (monomial) of $f$ with respect to $<_\Theta$ with the leading term (monomial) of the list $\Theta $ with respect to $<$.
	We said before the reason of this type of order of using with respect to $<$ and with respect to $<_\Theta$. The goal and the proof is similar to the goal and the proof of Lemma~\ref{TO give 1}.
	Note that the definition of the ordering $<_\Theta$ is in fact some generalization of the compatibility of the composition by $\Theta$  with the monomial ordering $<$, which is the key for the proof in both  results.
	
	\item In Lemma~\ref{MO6}, it will be stated that the composition operation commutes with the least common multiple computation if and only if every row of the exponent matrix $\Mat (LM_{<}(\Theta))$ has at most one non-zero element which will be indeed satisfied if the list $LM_{<}(\Theta)$ is permuted powering (equivalent to the compatibility of composition by $\Theta$
	with the non-divisibility with respect to $<$ as we know from Lemma~\ref{ND give 1}). As what we said for Lemma~\ref{MO5}, Lemma~\ref{MO6} is another strong form of Lemma~\ref{ND give 2} and plays the same role which this lemma played.
	The proof depends on simplifying the condition $$ \ \forall p \ \forall q [\lcm (p\circ (LM_{<}(\Theta)),q\circ (LM_{<}(\Theta)))  = \lcm (p,q) \circ (LM_{<}(\Theta))]$$
	to an equivalent linear-algebraic condition and proof that this new condition is equivalent to that every row of $\Mat (LM_{<}(\Theta))$ has at most one non-zero element. Note that the reason for not appearing of the ordering $<_\Theta$ here is that the only used leading monomial here is for the list $\Theta$ which defined using $<$ not $<_\Theta$ and the fact that $ LM_{<_\Theta}(p) = LM_{<}(p)=p$.
	
	\item In Lemma~\ref{MO7},  all the previous lemmas (Lemmas~\ref{MO4}, \ref{MO5} and \ref{MO6}) will be bundled into it under some restriction. It will be used frequently and essentially in proving Lemma \ref{MO8}, which forms the core of the proof of the main theorem. The choosing of the condition (list $LM_{<}(\Theta)$ is permuted powering) to be the condition (A) here is because this condition can play the rule played by all the (A)-conditions in the Lemmas~\ref{MO4}, \ref{MO5} and \ref{MO6} in the necessary part in these lemmas. This is because if the list $LM_{<}(\Theta)$ is permuted powering ($\equiv \det (\Mat (LM_{<}(\Theta)) \neq 0)$), then
	$<_\Theta$ is a monomial ordering and every row of $\Mat (LM_{<}(\Theta))$ has at most one non-zero element, so that all (B)-conditions of them are satisfied. For the other direction, if (B) here is true, then all (B)-conditions of these lemmas are true, which means that the matrix $\Mat (LM_{<}(\Theta))$ is non-singular with at most one non-zero element in each raw. This implies that $\Mat (LM_{<}(\Theta))$ must be a permuted diagonal matrix with exactly one non-zero element in each raw and column and thus the list $LM_{<}(\Theta)$ is permuted powering.
	
\end{itemize}

\begin{lemma}[\cite{Ho2}] \label{MO5}
	Fix a monomial ordering $<$ on the polynomial ring $K[x_1, \dots, x_n]$. Let $\Theta$ be a list of non-zero polynomials in $K[x_1, \dots, x_n]$. Let
	\begin{itemize}
		\item (A) The binary relation $<_\Theta$ is a monomial ordering;
		\item (B) For every $f \in K[x_1, \dots, x_n]$, we have
		\begin{enumerate}
			\item $LT_{<}(f\circ\Theta) = LT_{<_\Theta}(f)\circ (LT_{<}(\Theta))$.
			\item $LM_{<}(f\circ\Theta) = LM_{<_\Theta}(f)\circ (LM_{<}(\Theta))$.
		\end{enumerate}
	\end{itemize}
	Then (A) $\Rightarrow$ (B).
\end{lemma}
\begin{lemma}[\cite{Ho2}] \label{MO6}
	Fix a monomial ordering $<$ on the polynomial ring $K[x_1, \dots, x_n]$. Let $\Theta$ be a list of non-zero polynomials in $K[x_1, \dots, x_n]$, and  $p$ and $q$ are two monomials in $K[x_1, \dots, x_n]$. Let
	\begin{itemize}
		\item (A) Every row of the exponent matrix $\Mat (LM_{<}(\Theta))$ has at most one non-zero element; and
		\item (B) $ \ \forall p \ \forall q [\lcm (p\circ (LM_{<}(\Theta)),q\circ (LM_{<}(\Theta)))  = \lcm (p,q) \circ (LM_{<}(\Theta))]$.
	\end{itemize}
	Then (A) $\Longleftrightarrow$ (B).
\end{lemma}

\begin{lemma}[\cite{Ho2}] \label{MO7}
	Fix a monomial ordering $<$ on the polynomial ring $K[x_1, \dots, x_n]$. Let $\Theta$ be a list of non-zero polynomials in $K[x_1, \dots, x_n]$. Let $f$ be a polynomial and $p$ and $q$ be two monomials in $K[x_1, \dots, x_n]$. Let
	\begin{itemize}
		\item (A) The list $LM_{<}(\Theta)$ is permuted powering;
		\item (B) We have
		\begin{enumerate}
			\item The binary relation $<_\Theta$ is a monomial ordering.	
			\item $ \ \forall p \ \forall q [\lcm (p\circ (LM_{<}(\Theta)),q\circ (LM_{<}(\Theta)))  = \lcm (p,q) \circ (LM_{<}(\Theta))]$.	
			\item $LT_{<}(f\circ\Theta) = LT_{<_\Theta}(f)\circ (LT_{<}(\Theta))$.
			\item $LM_{<}(f\circ\Theta) = LM_{<_\Theta}(f)\circ (LM_{<}(\Theta))$.
		\end{enumerate}
	\end{itemize}
	Then (A) $\Longleftrightarrow$ (B).
\end{lemma}
As we commented before, Lemma~\ref{MO8} is the core for the proof of the main result of this section. The value of it is like the value of Lemma~\ref{not comp, suff} in the proof of the sufficiency part of the main result of the last two sections and the proofs of both are similar also. The compatibility of the composition with the monomial order exists here as a meaning of the definition of $<_\Theta$.
\begin{lemma}[\cite{Ho2}] \label{MO8} Fix a monomial ordering $<$ on the polynomial ring $K[x_1, \dots, x_n]$. Let $\Theta$ be a list of non-zero polynomials in $K[x_1, \dots, x_n]$. Let $G$ be a subset of non-zero polynomials in $K[x_1, \dots, x_n]$. Let
	\begin{itemize}
		\item (A) The list $LM_{<}(\Theta)$ is permuted powering;
		\item (B) $ \ \forall G \ [G. B._{<_\Theta}(G) \Rightarrow G. B._< (G\circ\Theta)]$.
	\end{itemize}
	Then (A) $\Rightarrow$ (B).
\end{lemma}
Finally, this is the main theorem in this section. As what did before, the only different between it and Lemma~\ref{MO8} is in the form of the condition which describes the problem studied in this section.
\begin{theorem}[\cite{Ho2}] \label{main result, different order} (Main Theorem for Different Monomial Ordering) Fix a monomial ordering $<$ on the polynomial ring $K[x_1, \dots, x_n]$. Let $\Theta$ be a list of non-zero polynomials in $K[x_1, \dots, x_n]$. Let $G$ and $F$ be two subsets of non-zero polynomials in $K[x_1, \dots, x_n]$. Let
	\begin{itemize}
		\item (A) The list $LM_{<}(\Theta)$ is permuted powering;
		\item (B) $ \ \forall F\ \forall G \ [G. B._{<_\Theta}(G,\langle F \rangle) \Rightarrow G. B._< (G\circ\Theta, \langle F\circ\Theta \rangle)]$.
	\end{itemize}
	Then (A) $\Rightarrow$ (B).
\end{theorem}
We end with some examples given also by Hong at the end of his paper \cite{Ho2}.

\begin{example}[\cite{Ho2}] \label{ex1,2}
	Every composition such that $LM_{<}(\Theta)$ is permuted powering is a solution of the studied problem in this section. This means that Examples~\ref{ex1,1}, \ref{ex2,1} and \ref{ex3,1} given at the end of Section~\ref{Gr\"obner Bases Under Composition2} can be considered as  examples of this section.
\end{example}
Now, since the answer of our problem is not complete and it only given a sufficient condition for this type of commutativity, it can be found some other examples with some special conditions. Hong gave some of these examples to explain why his results are important. We put some parts of them without his calculations and proofs.

\begin{example}[\cite{Ho2}]
	Fix any monomial ordering $<$. Let $F$ be a finite set of monomials and let $\Theta $ be such that $LM_{<}(\Theta)$ is permuted powering. Then $F\circ\Theta$ is a Gr\"obner basis with respect to $<$. This follows immediately from the main theorem and the fact that $F$ is already a Gr\"obner basis of the ideal $\langle F \rangle$ with respect to $<$.
\end{example}

\begin{example}[\cite{Ho2}]
	This example is with the special case $ < = <_\Theta$. Let  $ < $ be the graded reverse lexicographic ordering with $x_3 < x_2 < x_1$, and $$\Theta=( (x_1^2 + x_2 +x_3)^3, (x_1 + x_2^2 +x_3)^3, (x_1 + x_2 +x_3^2)^3).$$ In general, whenever $\Mat (LM_{<}(\Theta))$ is a diagonal matrix with the same diagonal entries, we have that $ < \, = \, <_\Theta$.
\end{example}
\begin{example}[\cite{Ho2}] \label{ex4,2}
	Now an example such that $ < \neq <_\Theta$. Let  $ < $ be the weighted reverse lexicographic ordering defined as $$x_1^{a_1} x_2^{a_2}x_3^{a_3} < x_1^{b_1} x_2^{b_2}x_3^{b_3} \Longleftrightarrow 2a_1+3a_2+6a_3< 2b_1+3b_2+6b_3,$$ where the tie will be broken by the reverse lexicographic ordering where $x_3 < x_2 < x_1$, and $$\Theta=( (x_1^3 + x_2 +x_3)^3, (x_1 + x_2^2 +x_3)^3, (x_1 + x_2 +x_3)^3).$$ The calculations will give us that the ordering $<_\Theta$, which will be a monomial ordering since $$LM_{<}(\Theta)= (x_1^9,x_2^6,x_3^3)$$ is permuted powering, is the graded reverse lexicographic ordering. Clearly, the main theorem can be used here.
\end{example}

\section{Reduced Gr\"obner Bases under Composition}\label{Reduced Gr\"obner Bases Under Composition}
The second open question given by Hong at the end of his first paper \cite{Ho1} in our subject was:

\emph{When does a composition commute with the reduced Groebner bases computation?}

In other words, which conditions must we put on the list $\Theta$ to make every reduced Gr\"obner basis $G$ with respect to some monomial $<$ for an ideal $I$ composes using $\Theta$ to a reduced Gr\"obner basis $G\circ\Theta$ with respect to the same $<$ for the ideal $\langle I\circ\Theta \rangle$.
Hong put this as an open question that is not solved using his main result proved in \cite{Ho1} because of:

\emph{One can easily construct an example that shows that the two conditions given in
	this paper are not sufficient.}

He commented:

\emph{ An answer to this question will shed new light on the notion of `reduced'.}

In 1998, J. Guti\'errez and R. Rubio San Miguel in  \cite{GR} gave a complete answer to this question. They proved that for every reduced Gr\"obner basis $G$ with respect to $<$, $G \circ\Theta$ is a reduced Gr\"obner basis with respect to the same monomial ordering $<$ if and only if the composition by $\Theta$ is compatible with the monomial ordering $<$ and $\Theta$ is a list of permuted univariate and monic polynomials. This result is different of Hong's result by the additional two conditions: univariate and monic. Every new condition will be need for a reason comes from the meaning of the word `reduced'. They also studied the third open question given and solved by Hong for reduced Gr\"obner bases and gave a sufficient condition to determine when composition commutes with reduced Gr\"obner bases computation under some possibly different monomial orderings where the second possibly different monomial ordering is which the ordering defined in Definition~\ref{Composition on Ordering}. Note that, as we said before in Section~\ref{Gr\"obner Bases Under Composition3}, this is a generalization in some way for the first question.  In addition to these studies for reduced Gr\"obner bases case, they studied these two questions for minimal Gr\"obner bases case and proved some results near to the results proved by Hong in his two papers \cite{Ho1,Ho2}.
In this section, we will study this paper and describe their work with giving our comments on it. We will start by giving the mathematical notations of the studied problems in this section. The reader will notice that there is no differences in the way of describing the problems in this chapter because all authors (included us) which studied this subject followed Hong's method of writing.
\begin{definition}[\cite{GR}]
	(Commutativity with Composition of Minimal and Reduced Gr\"obner Basis) Fix a monomial ordering $<$ on the polynomial ring $K[x_1, \dots, x_n]$. We say that the composition by $\Theta $ commutes with the computation of minimal Gr\"obner bases with respect to $<$ if the following formula is true for $\Theta $:
$$ \ \forall I \ \forall G \ [M. G. B._< (G, I) \Rightarrow M. G. B._< (G\circ\Theta, \langle I\circ\Theta \rangle)].$$
	We say that the composition by $\Theta $ commutes with the computation of reduced Gr\"obner bases with respect to $<$ if the following formula is true for $\Theta $:
 $$ \ \forall I \ \forall G \ [R. G. B._< (G, I) \Rightarrow R. G. B._< (G\circ\Theta, \langle I\circ\Theta \rangle)].$$
	Also, we say that the composition by $\Theta $ commutes with the computation of reduced sets with respect to $<$ if the following formula is true for $\Theta $:
 $$ \ \forall F \ [R._< (F) \Rightarrow R._< (F\circ\Theta)].$$
\end{definition}
The reader must notice and differ later between the places of using of reduced sets and the places of using reduced Gr\"obner bases since not every reduced set need to be a reduced Gr\"obner basis. Using reduced sets in some results is to give a strong or a two sided result which may not hold for reduced Gr\"obner bases. The meaning of the additional condition `univariate' will be given in the next definition. This will be used later in the proofs to solve the problems in computations which come from the condition:
$$\text{For all} ~f \in F,~ \text{no monomial of}~ f ~\text{lies in} ~ \langle LT_{<}(F - \{f\})  \rangle,$$ which is the main condition in the definition of reduced sets. Using  the condition `monic' later is for
$$LC_{<}(f) = 1~ \text{for all} ~f \in F.$$ This is very important and we will see later that without these two conditions, the list $\Theta$ will only satisfy the conditions of Hong's result and this means that the reduced Gr\"obner basis will compose to a Gr\"obner basis which need not be reduced. For the minimal Gr\"obner bases cases, we will see that the only  additional condition which will be needed is the condition `monic'.
\begin{definition}[\cite{GR}]
	(Permuted Univariate and Monic Polynomials) Fix a monomial ordering $<$ on the polynomial ring $K[x_1, \dots, x_n]$. Let $\Theta$ be a list of non-zero polynomials in $K[x_1, \dots, x_n]$. We say that $\Theta$ is a list of permuted univariate and monic polynomials if $$\Theta= (f_1(x_{\pi_1}),\dots,f_n(x_{\pi_n})),$$ where $LT_{<}(f_i) =x^{\lambda_i}_{\pi_i}$ with $\pi_j = \pi(j)$ for $j \in \{1,\dots, n\}$, and $\pi$  is a permutation of $(1,\dots, n)$ and $\lambda_1,\dots, \lambda_n > 0$.
\end{definition}
 Theorem~\ref{Reduced Gr\"obner Bases Under Composition, main result} is the main theorem of commutativity of reduced Gr\"obner bases with respect to $<$ under Composition. It set that, as we said before, the composition by $\Theta $ commutes with the computation of reduced Gr\"obner basis with respect to $<$ if and only if the composition by $\Theta$ is compatible with the monomial ordering $<$ and $\Theta$ is a list of permuted univariate and monic polynomials. The proof of this theorem needs the following three results:
\begin{itemize}
	\item Lemma~\ref{Red set give nondi} is an easy proved result that says that if the composition by $\Theta $ commutes with the computation of reduced sets with respect to $<$, then it is compatible with the non-divisibility with respect to $<$. It will be used in the proof of Proposition~\ref{Red set equa pum} and Lemma~\ref{Red set give To} with the fact that the compatibility with the non-divisibility with respect to $<$ is equivalent to that $LM_{<}(\Theta)$ is permuted powering. Note that the result will be true for reduced Gr\"obner bases.
	
	\item Proposition~\ref{Red set equa pum} itself is the answer if the question is for the conditions on $\Theta$ which make the composition by it commutes with the computation of the reduced sets. This will  happen if and only if $\Theta$ is a list of permuted univariate and monic polynomials. The use of Lemma~\ref{Red set give nondi} is for the second direction to make the composition preserves the non-divisibility. Using reduced sets here is very important to give a strong two sided result since this proposition is valued only for the first direction for reduced Gr\"obner bases as we will say later.
	
	\item  Lemma~\ref{Red set give To}, the first use of reduced Gr\"obner bases where the commutativity here implies the compatibility with the monomial ordering. The proof here is directly using some calculations like what Hong used in his proofs but with different way of defining the used sets beside the result proved in Lemma~\ref{Red set give nondi}. This result with Hong's result given in Theorem~\ref{Gr\"obner Bases Under Composition, main result} improve Proposition~\ref{Red set equa pum} to give us Theorem~\ref{Reduced Gr\"obner Bases Under Composition, main result} for reduced Gr\"obner bases case.
\end{itemize}
Now, the proof of sufficiency part of Theorem~\ref{Reduced Gr\"obner Bases Under Composition, main result} will be using of Theorem~\ref{Gr\"obner Bases Under Composition, main result} and Proposition~\ref{Red set equa pum}. Let $G$ be a reduced Gr\"obner basis with respect to some $<$ for the ideal $I$. That means that $G$ is also a reduced set,
and thus by Proposition~\ref{Red set equa pum}, $G\circ\Theta$ is also a reduced set. Now the composition by $\Theta$ satisfied the conditions of
Theorem~\ref{Reduced Gr\"obner Bases Under Composition, main result}, so that $G\circ\Theta$ is a Gr\"obner basis with respect to $<$ for the ideal $ \langle I\circ\Theta \rangle$. Hence $G\circ\Theta$ is a reduced Gr\"obner basis with respect to $<$ for the ideal $\langle I\circ\Theta \rangle$.
For the proof of necessity, one half of it will need Proposition~\ref{Red set equa pum} and the other half is by using Lemma~\ref{Red set give To}.
Although Proposition~\ref{Red set equa pum} gives us the equivalently between the commutativity of composition by $\Theta$ with the computation of reduced sets
and that $\Theta$ is a list of permuted univariate and monic polynomials, the proof of the first direction of it (which we only need here) is the same
if we assume that the commutativity holds for reduced Gr\"obner bases only since the set used in the proof is also a reduced Gr\"obner bases.
Lemma~\ref{Red set give To} is clear and shows us the second half directly.

\begin{lemma}[\cite{GR}] \label{Red set give nondi}
	Fix a monomial ordering $<$ on the polynomial ring $K[x_1, \dots, x_n]$. Let $\Theta$ be a list of non-zero polynomials in $K[x_1, \dots, x_n]$. Let
	\begin{itemize}
		\item (A) Composition by $\Theta $ commutes with the computation of reduced sets with respect to $<$;
		\item (B) Composition by $\Theta $ is compatible with the non-divisibility with respect to $<$.
	\end{itemize}
	Then (A) $\Rightarrow$ (B).	
\end{lemma}
\begin{proposition}[\cite{GR}] \label{Red set equa pum}
	Fix a monomial ordering $<$ on the polynomial ring $K[x_1, \dots, x_n]$. Let $\Theta$ be a list of non-zero polynomials in $K[x_1, \dots, x_n]$. Let
	\begin{itemize}
		\item (A) Composition by $\Theta $ commutes with the computation of reduced sets with respect to $<$;
		\item (B) $\Theta$ is a list of permuted univariate and monic polynomials.
	\end{itemize}
	Then (A) $\Longleftrightarrow$ (B).	
\end{proposition}
\begin{lemma}[\cite{GR}] \label{Red set give To}
	Fix a monomial ordering $<$ on the polynomial ring $K[x_1, \dots, x_n]$. Let $\Theta$ be a list of non-zero polynomials in $K[x_1, \dots, x_n]$. Let
	\begin{itemize}
		\item (A) Composition by $\Theta $ commutes with the computation of reduced Gr\"obner bases with respect to $<$;
		\item (B) Composition by $\Theta $ is compatible with the monomial ordering $<$.
	\end{itemize}
	Then (A) $\Rightarrow$ (B).	
\end{lemma}
\begin{theorem}[\cite{GR}] \label{Reduced Gr\"obner Bases Under Composition, main result}
	(Main theorem of commutativity of reduced Gr\"obner bases with respect to $<$ under Composition) Fix a monomial ordering $<$ on the polynomial ring $K[x_1, \dots, x_n]$. Let $\Theta$ be a list of non-zero polynomials in $K[x_1, \dots, x_n]$, then the following are equivalent
	\begin{itemize}
		\item (A) Composition by $\Theta $ commutes with the computation of reduced Gr\"obner basis with respect to $<$;
		\item (B) Composition by $\Theta $ is
		\begin{enumerate}
			\item compatible with the monomial ordering $<$; and
			\item $\Theta$ is a list of permuted univariate and monic polynomials.
		\end{enumerate}	
	\end{itemize}	
\end{theorem}
Repeating what we said before, it is important to point out here that Theorem~\ref{Reduced Gr\"obner Bases Under Composition, main result} and
Theorem~\ref{Gr\"obner Bases Under Composition, main result} are intrinsically different. However, for minimal Gr\"obner bases it can easily
obtain (using Theorem~\ref{Reduced Gr\"obner Bases Under Composition, main result})
that a minimal Gr\"obner basis with respect to $<$ computation commutes with composition if and only if the composition is compatible
with the monomial ordering $<$ and the non-divisibility and $LC_{<}(\Theta)= (1, \dots,1)$.
For the first part, that is moving from the commutativity of Gr\"obner bases to commutativity of reduced Gr\"obner bases, we need compatibility with the non-divisibility
(which  clearly holds) besides that $LC_{<}(\Theta)= (1, \dots,1)$. The reverse direction will be using the fact that from any Gr\"obner basis $G$,
we can remove elements from it then dividing on the leading coefficients to get some minimal Gr\"obner basis $G'$ with respect to the same monomial ordering
for the ideal $\langle G \rangle$. $LC_{<}(\Theta)= (1, \dots,1)$ will be proved using that $\{x_i\}$ is a minimal Gr\"obner basis for any $i$.
\begin{theorem}[\cite{GR}]
	(Main theorem of commutativity of minimal Gr\"obner bases with respect to $<$ under Composition) Fix a monomial ordering $<$ on the polynomial ring $K[x_1, \dots, x_n]$. Let $\Theta$ be a list of non-zero polynomials in $K[x_1, \dots, x_n]$, then the following are equivalent
	\begin{itemize}
		\item (A) Composition by $\Theta $ commutes with the computation of Gr\"obner bases with respect to  $<$ and $LC_{<}(\Theta)= (1, \dots,1)$;
		\item (B) Composition by $\Theta $ commutes with the computation of minimal Gr\"obner bases with respect to $<$.
	\end{itemize}	
\end{theorem}
We will study now the behaviour
of minimal and reduced Gr\"obner bases under composition of polynomials (possibly) under different monomial orderings. Remember that from Definition~\ref{Composition on Ordering} and Lemma~\ref{MO7}, the binary relation $<_\Theta$  defined by $$\ \forall p \ \forall q \ [~p ~ <_\Theta q \Longleftrightarrow p\circ (LM_{<}(\Theta)) < q\circ (LM_{<}(\Theta))]$$
is a monomial ordering if the list $LM_{<}(\Theta)$ is permuted powering. Note that when $\Theta$ is compatible with the monomial ordering $<$, the binary relation $<_\Theta$ is exactly $<$. The reader will notice that, besides the condition $LC_{<}(\Theta)= (1, \dots,1)$, there is no additional condition on the compatibility of the composition with minimal Gr\"obner bases under these monomial orderings. For the reduced Gr\"obner bases case, we need only the condition that $\Theta$ is a list of permuted univariate and monic polynomials since the compatibility with the monomial orderings $<$ and $<_\Theta$ is hold from the definition of $<_\Theta$. The following two theorems describe what we said in this paragraph.
\begin{theorem}[\cite{GR}] Fix a monomial ordering $<$ on the polynomial ring $K[x_1, \dots, x_n]$. Let $\Theta$ be a list of non-zero polynomials in $K[x_1, \dots, x_n]$. Let $G$ and $I$ be two subsets of non-zero polynomials in $K[x_1, \dots, x_n]$. Let
	\begin{itemize}
		\item (A) The list $LM_{<}(\Theta)$ is permuted powering and $LC_{<}(\Theta)= (1, \dots,1)$;
		\item (B) $ \ \forall I\ \forall G \ [M. G. B._{<_\Theta}(G,I) \Rightarrow M. G. B._< (G\circ\Theta, \langle I\circ\Theta \rangle)]$.
	\end{itemize}
	Then (A) $\Rightarrow$ (B).
\end{theorem}

\begin{theorem}[\cite{GR}] Fix a monomial ordering $<$ on the polynomial ring $K[x_1, \dots, x_n]$. Let $\Theta$ be a list of non-zero polynomials in $K[x_1, \dots, x_n]$. Let $G$ and $I$ be two subsets of non-zero polynomials in $K[x_1, \dots, x_n]$. Let
	\begin{itemize}
		\item (A) $\Theta$ is a list of permuted univariate and monic polynomials;
		\item (B) $ \ \forall I\ \forall G \ [R. G. B._{<_\Theta}(G,I) \Rightarrow R. G. B._< (G\circ\Theta, \langle I\circ\Theta \rangle)]$.
	\end{itemize}
	Then (A) $\Rightarrow$ (B).
\end{theorem}
We give some examples (see \cite{GR}) for the results giving in this section.

\begin{example}[\cite{GR}] \label{ex1,3}
	Let $<$ be the lexicographic ordering. Then, every composition of the form $$\Theta= (\theta_1(x_1),\dots,\theta_n(x_n)),$$ with $LM_<(\theta_i)=x_i^{\lambda_i}$ and $0 < \lambda_i $ for every $i$, commutes with the computation of reduced Gr\"obner bases with respect to $<$.
\end{example}

\begin{example}[\cite{GR}] \label{ex2,3}
	Let $<$ be the graded lexicographic ordering. Then, every composition of the form $$\Theta= (\theta_1(x_1),\dots,\theta_n(x_n)),$$ with $LM_<(\theta_i)=x_i^{\lambda}$ for every $i$ and $0 < \lambda$, commutes with the computation of reduced Gr\"obner bases with respect to $<$.	
\end{example}

\begin{example}[\cite{GR}] \label{ex3,3}
	Fix any monomial ordering $<$ and some positive real numbers $\alpha_1, \dots, \alpha_n$. Consider the monomial ordering $<_L$ defined as follows. Let $p= x_1^{a_1} \dots x_n^{a_n}$	and $q= x_1^{b_1} \dots x_n^{b_n}$, then
	
	\[p <_L q :\Longleftrightarrow
\begin{cases}
	\alpha_1a_1+\dots+\alpha_na_n < \alpha_1b_1+\dots+\alpha_nb_n, \\
	\text{or} \\
	\alpha_1a_1+\dots+\alpha_na_n = \alpha_1b_1+\dots+\alpha_nb_n,~\text{and}~ p < q.\\
\end{cases}\]
	If $\Theta$ is compatible with the monomial ordering $<$, then every composition of the form $$\Theta= (\theta_1(x_1),\dots,\theta_n(x_n)),$$ with $LM_{<_L}(\theta_i)=x_i^{\lambda}$ for every $i$ and $0 < \lambda$, commutes with the computation of reduced Gr\"obner bases with respect to $<_L$. In fact, Example~\ref{ex2,3} is a particular case by taking $a_i =1$ for any $i$ and $<$ the lexicographic ordering.
\end{example}

\begin{example}[\cite{GR}] \label{ex4,3}
	Remember the monomial ordering $<$ defined in Example~\ref{ex3,1} by $$ p < q \Longleftrightarrow a_1 +\sqrt{2} a_2 < b_1 +\sqrt{2} b_2,$$ where
	$p=x^{a_1}y^{a_2}$ and $q=x^{b_1}y^{b_2}$ are two monomials in $K[x,y]$.  Then every composition of the form $$\Theta= (\theta_1(x),\theta_2(y)),$$ where $LM_{<}(\Theta)= (x^\lambda,y^\lambda) $ with $0 < \lambda$, commutes with the computation of reduced Gr\"obner bases with respect to $<$.
	
\end{example}

\begin{example}[\cite{GR}] \label{ex5,3}
	With the same monomial ordering $<$ defined above in Example~\ref{ex4,3}, every composition of the form $$\Theta= (\theta_1(y),\theta_2(x)),$$ where $LM_{<}(\Theta)= (y^\lambda,x^{2\lambda}) $ with $0 < \lambda$, commutes with the computation of reduced Gr\"obner bases with respect to $<$. This case involves a non-trivial permutation of variables where all previous cases involve trivial permutations of variables.
\end{example}
\section{Gr\"obner Bases under Composition: Different Polynomial Rings and Monomial Orderings}\label{Gr\"obner Bases Under Composition: Different Polynomial Rings and Monomial Orderings}
\sectionmark{Different Polynomial Rings and Monomial Orderings}
In Sections~\ref{Gr\"obner Bases Under Composition1} and \ref{Gr\"obner Bases Under Composition2} we studied Hong's paper \cite{Ho1} with the following problem: when does composition commute with the computation of Gr\"obner bases under some monomial ordering? Hong's second paper \cite{Ho2} was studied in Section~\ref{Gr\"obner Bases Under Composition3} with some other look of this problem by change the monomial ordering in one side under some special definition.
In Section~\ref{Reduced Gr\"obner Bases Under Composition}, we studied the paper \cite{GR} of J. Guti\'errez and R. Rubio San Miguel with similar problems for minimal and reduced Gr\"obner bases. In this section, we will study the work of Z. Liu and M. Wang \cite{LZW}  published in 2001. In that paper, they studied the same previous problems but with some generalizations of the definitions of commutativity and compatibility given in the beginnings of these previous sections (Section~\ref{Gr\"obner Bases Under Composition1}--Section~\ref{Reduced Gr\"obner Bases Under Composition}). These generalizations are in two ways:
\begin{itemize}
	\item The first is by using two polynomial rings in the definitions, the first is $K[x_1, \dots, x_n]$ for polynomials before the composition and the other is $K[y_1, \dots, y_m]$ to be used after it. To make this make sense, we must take $n \leq m$ and $\Theta = (\theta_1, \dots, \theta_n)$, where $\theta_i \in K[y_1, \dots, y_m]$ for any $i$. Note that in all previous sections we work over the case $n = m$.
	
	\item The second way is just what we said before in Section~\ref{Gr\"obner Bases Under Composition3} about the second general meaning of the words `another monomial ordering' appeared in the question:
	
	\emph{Let $G$ be a Gr\"obner basis for $I$ with respect to $<$. When is $G \circ \Theta$ is a Gr\"obner basis for $\langle I \circ \Theta \rangle$ (possibly with respect to another monomial  ordering $<'$)?}
	
	\noindent given by Hong in \cite{Ho1}. We commented by
	
	\emph{When both of the two monomial orderings are arbitrary without any special conditions on any of them. In other words we choose two monomial orderings from the beginning and use one to define the Gr\"obner basis $G$ with respect to it before the composition and use the other ordering for $G \circ \Theta$ after this operation.}
	
	 Remember that what was studied in Sections~\ref{Gr\"obner Bases Under Composition3} and \ref{Reduced Gr\"obner Bases Under Composition} is a special case with some special definition of one of these monomial orderings for the general, minimal and reduced Gr\"obner bases over the case $n = m$ also.
\end{itemize}	
These ways were studied in this paper together, that is the generalized definitions and results given in it is by using both of previous ways in the same time as we will start to see in the two definitions given after this paragraph. Z. Liu and M. Wang proved first in Theorem~\ref{Gr\"obner Bases Under Composition,diff monomial ordering, main result}, that over these two ways of generalizations (different polynomial rings $K[x_1, \dots, x_n]$ and $K[y_1, \dots, y_m]$ and different monomial orderings $<_1$ and $<_2$ before and after the composition), the composition by $\Theta$ commutes with the computation of Gr\"obner bases with respect to $<_1$ and $<_2$ if and only if this composition by is compatible with $<_1$ and $<_2$ and for all $i\neq j$, the monomials $LM_{<_2}(\theta_i)$ and $LM_{<_2} (\theta_j)$ are relatively prime. This thing is more complicated for the reduced Gr\"obner bases case, but they gave some good answers using different sufficient and necessary condition for this case of commutativity as we will see in Theorem~\ref{reduced Gr\"obner Bases Under Composition,diff monomial ordering, main result}, Theorem~\ref{reduced Gr\"obner Bases Under Composition,diff monomial ordering, main result2}, and Corollary~\ref{reduced Gr\"obner Bases Under Composition,diff monomial ordering, coro}. As done in every section, the needed definitions will be set firstly in the following two definitions.
\begin{definition}[\cite{LZW}]
	Fix a monomial ordering $<_1 $ on the polynomial ring $K[x_1, \dots, x_n]$ and a monomial ordering $<_2$ on the polynomial ring $K[y_1, \dots, y_m]$. Let $\Theta = (\theta_1, \dots, \theta_n)$ be a list of $n$ non-zero polynomials in $K[y_1, \dots, y_m]$.
	\begin{enumerate}
		\item We say that the composition by $\Theta $ commutes with the computation of Gr\"obner bases with respect to $<_1$ and $<_2$ if the following formula is true for $\Theta $: $$ \ \forall I \ \forall G \ [G. B._{<_1} (G, I) \Rightarrow G. B._{<_2} (G\circ\Theta, \langle I\circ\Theta \rangle)].$$
		\item Also, we say that the composition by $\Theta $ commutes with the computation of reduced Gr\"obner bases with respect to $<_1$ and $<_2$ if the following formula is true for $\Theta $: $$ \ \forall I \ \forall G \ [R. G. B._{<_1} (G, I) \Rightarrow R. G. B._{<_2} (G\circ\Theta, \langle I\circ\Theta \rangle)].$$
	\end{enumerate}
\end{definition}
\begin{definition}[\cite{LZW}]
	Fix a monomial ordering $<_1 $ on the polynomial ring $K[x_1, \dots, x_n]$ and a monomial ordering $<_2$ on the polynomial ring $K[y_1, \dots, y_m]$. Let $\Theta = (\theta_1, \dots, \theta_n)$ be a list of $n$ non-zero polynomials in $K[y_1, \dots, y_m]$.
	\begin{enumerate}
		\item We say that the composition by $\Theta $ is compatible with the monomial orderings $<_1$ and  $<_2$ if for all monomials $p$ and $q$ in $K[x_1, \dots, x_n]$, the following formula is true for $\Theta $: $$ \ \forall p \ \forall q \ [p <_1 q \Rightarrow p\circ (LM_{<_2}(\Theta)) <_2 q\circ (LM_{<_2}(\Theta))].$$
		\item We say that the composition by $\Theta $ is compatible with the non-divisibility with respect to $<_1$ and  $<_2$ if for all monomials $p$ and $q$ in $K[x_1, \dots, x_n]$, the following formula is true for $\Theta $: $$ \ \forall p \ \forall q \ [p\nmid q \Rightarrow p\circ (LM_{<_2}(\Theta)) \nmid q\circ (LM_{<_2}(\Theta))].$$
	\end{enumerate}
\end{definition}
The following lemma seems to be identical to Lemma ~\ref{propert 1} and actually  both of them are used in the same way and their proofs are quite similar. Nevertheless, it is important to prove it, since here we have different polynomial rings and monomial orderings on the two sides. This lemma shows us how to compute the leading monomials and terms after the composition and it is used in the proof of Lemma~\ref{diff TO give 1}.

\begin{lemma} [\cite{LZW}]	\label{prop1}
	Fix a monomial ordering $<_1 $ on the polynomial ring $K[x_1, \dots, x_n]$ and a monomial ordering $<_2$ on the polynomial ring $K[y_1, \dots, y_m]$.
	Let $\Theta = (\theta_1, \dots, \theta_n)$ be a list of $n$ non-zero polynomials in $K[y_1, \dots, y_m]$. Then for any monomial $p$ in $K[x_1, \dots, x_n]$ we have
	\begin{enumerate}
		\item $LT_{<_2}(p\circ\Theta) = p\circ (LT_{<_2}(\Theta))$.
		\item $LM_{<_2}(p\circ\Theta) = p\circ (LM_{<_2}(\Theta))$.
	\end{enumerate} 	
\end{lemma}
The same comment before Lemma~\ref{prop1} can be told here because Lemma~\ref{diff TO give 1} is similar to several previous lemmas like Lemma~\ref{TO give 1} and Lemma~\ref{MO7}.
The following lemma is a strong generalisation of these previous lemmas, and this not only because as we said before we have different polynomial rings and monomial orderings on the two sides.
This is because this lemma is also a two sided result and the proof of the other side can be used to generalise Lemma~\ref{TO give 1} and Lemma~\ref{MO7} to be two sided. By these two sided results, we can get a method to decide if the composition by $\Theta $ is compatible with the monomial orderings or not using calculations and then know which forms of $\Theta $ satisfy this condition.
\begin{lemma} [\cite{LZW}] \label{diff TO give 1}
	Fix a monomial ordering $<_1 $ on the polynomial ring $K[x_1, \dots, x_n]$ and a monomial ordering $<_2$ on the polynomial ring $K[y_1, \dots, y_m]$. Let $\Theta = (\theta_1, \dots, \theta_n)$ be a list of $n$ non-zero polynomials in $K[y_1, \dots, y_m]$. Let
	\begin{itemize}
		\item (A) Composition by $\Theta $ is compatible with the monomial orderings $<_1$ and $<_2$;
		\item (B) For every non-zero $f \in K[x_1, \dots, x_n]$, we have
		\begin{enumerate}
			\item $LT_{<_2}(f\circ\Theta) = LT_{<_1}(f)\circ (LT_{<_2}(\Theta))$.
			\item $LM_{<_2}(f\circ\Theta) = LM_{<_1}(f)\circ (LM_{<_2}(\Theta))$.
			\item If $LC_{<_2}(\Theta)= (1, \dots,1)$ and $LC_{<_1}(f)=1$, then $LC_{<_2}(f\circ\Theta) = 1$.
		\end{enumerate}
	\end{itemize}
	Then (A)$\Longleftrightarrow$ (B).
\end{lemma}
The new thing in the next lemma (besides the different polynomial rings and monomial orderings), which in not found in Lemma~\ref{ND give 1}, is the existing of the monomials $t_i$ from the polynomial ring $K [\{y_1, \dots, y_m\} - \{y_{\pi_1},\dots,y_{\pi_n}\}]$ in the form of the list $LM_{<_2}(\Theta)$ because we work over the general case $n \leq m$	and this will give us some free variables not used after finishing the variables used in the injective map $\pi$. The proof of both lemmas is similar.

\begin{lemma}[\cite{LZW}]
	Fix a monomial ordering $<_1 $ on the polynomial ring $K[x_1, \dots, x_n]$ and a monomial ordering $<_2$ on the polynomial ring $K[y_1, \dots, y_m]$. Let $\Theta = (\theta_1, \dots, \theta_n)$ be a list of $n$ non-zero polynomials in $K[y_1, \dots, y_m]$. Let
	\begin{itemize}
		\item (A) Composition by $\Theta $ is compatible with the non-divisibility with respect to $<_1$ and $<_2$;
		\item (B) the list $LM_{<_2}(\Theta)= (y^{\lambda_1}_{\pi_1}\cdot t_1,\dots,y^{\lambda_n}_{\pi_n}\cdot t_n)$ for some injective map $\pi$ from $(1,\dots, n)$ to $(1,\dots, m)$, and for each $i$, $t_i$ is a monomial in the variables $ \{y_1, \dots, y_m\} - \{y_{\pi_1},\dots,y_{\pi_n}\}$, where $\lambda_1,\dots, \lambda_n > 0$.
	\end{itemize}
	Then (A)$\Longleftrightarrow$ (B).
\end{lemma}
Since Lemma~\ref{diff TO give 1} is a strong generalisation of Lemma~\ref{TO give 1} and Lemma~\ref{MO7},
the following lemmas are a very strong two sided generalisation of Lemma~\ref{ND give 2} and Lemma~\ref{MO7}
for the same reasons and also because we have a new equivalent condition here used to simplify the proof and gives us the method to generalise these lemmas.
This condition shows us how we can compute the greatest common divisor of the leading monomials after the composition.
Also, it gives us the relatively prime condition of the monomials $LM_{<_2}(\theta_i)$ and $LM_{<_2} (\theta_j)$ for all $i\neq j$ directly (which is the key for the permuted powering condition as we know) without working with a helping pre-lemmas or using long proofs over several cases.
\begin{lemma}[\cite{LZW}] \label{diff ND give 2}
	Fix a monomial ordering $<_1 $ on the polynomial ring $K[x_1, \dots, x_n]$ and a monomial ordering $<_2$ on the polynomial ring $K[y_1, \dots, y_m]$. Let $\Theta = (\theta_1, \dots, \theta_n)$ be a list of $n$ non-zero polynomials in $K[y_1, \dots, y_m]$. Let
	\begin{itemize}
		\item (A) for all $i\neq j$, the monomials $LM_{<_2}(\theta_i)$ and $LM_{<_2} (\theta_j)$ are relatively prime;
		\item (B) $ \ \forall p \ \forall q  \ $ $ [\lcm (p\circ (LM_{<_2}(\Theta)),q\circ (LM_{<_2}(\Theta)))  = \lcm (p,q) \circ (LM_{<_2}(\Theta))]$.
		\item (C) $ \ \forall p \ \forall q  \ $ $ [\gcd (p\circ (LM_{<_2}(\Theta)),q\circ (LM_{<_2}(\Theta)))  = \gcd (p,q) \circ (LM_{<_2}(\Theta))]$.
	\end{itemize}
	Then (A) $\Longleftrightarrow$ (B) $\Longleftrightarrow$ (C).
\end{lemma}
This is the last step before setting the first main result of this section. One direction of the following lemma is just the contrapositive of Lemma~\ref{not gb1} given by Hong but in a general form by working with both of the two polynomials $f$ and $g$. The other direction is the converse of that contrapositive. The proof depends on the greatest common divisor condition given in Lemma~\ref{diff ND give 2} above. Also, Lemma~\ref{gcd} which will be given in the next section, will be used here and the authors of the studied paper \cite{LZW} in this section put it as a claim in the proof.
\begin{lemma}[\cite{LZW}]
	Fix a monomial ordering $<_2$ on $K[y_1, \dots, y_m]$. Let $f$ and $g$ be  two non-zero polynomials in $K[y_1, \dots, y_m]$. Then $ \{f,g\}$, $ \{f+1,g\}$ and $ \{f,g+1\}$ are Gr\"obner bases with respect to $<_2$ if and only if
	\[\gcd (LM_{<}(f), LM_{<}(g))=  1.\]
\end{lemma}
The following result is the first main theorem in this section. It will be the end of the story of the commutativity of the composition with the computation of Gr\"obner bases in the general case and gives us a generalization of all main results given by Hong in his two papers \cite{Ho1,Ho2} because it will work for any polynomial ring and monomial ordering in any side before and after the composition. The forms of all these main results are the same, but some time some condition is used and other time an equivalent to it is used as done here with the condition: for all $i\neq j$, the monomials $LM_{<_2}(\theta_i)$ and $LM_{<_2} (\theta_j)$ are relatively prime.
The results of Hong can be obtained here as  follows:
\begin{itemize}
	\item If $m=n$ and $<_1$ is equal to $<_2$, then Theorem~\ref{Gr\"obner Bases Under Composition,diff monomial ordering, main result} is the same as Theorem~\ref{Gr\"obner Bases Under Composition, main result}.
	\item If $m=n$ and $<_1$ is equal to $<_\Theta$, where $<_2$ is written as $<$, then Theorem~\ref{Gr\"obner Bases Under Composition,diff monomial ordering, main result} is a strong two sided generalization of Theorem~\ref{main result, different order}.
\end{itemize}
\begin{theorem}[\cite{LZW}] \label{Gr\"obner Bases Under Composition,diff monomial ordering, main result}
	Fix a monomial ordering $<_1 $ on the polynomial ring $K[x_1, \dots, x_n]$ and a monomial ordering $<_2$ on the polynomial ring $K[y_1, \dots, y_m]$. Let $\Theta = (\theta_1, \dots, \theta_n)$ be a list of $n$ non-zero polynomials in $K[y_1, \dots, y_m]$. Then the following are equivalent
	\begin{itemize}
		\item (A) Composition by $\Theta $ commutes with the computation of Gr\"obner bases with respect to $<_1$ and $<_2$;
		\item (B) Composition by $\Theta $ is
		\begin{enumerate}
			\item compatible with the monomial orderings $<_1$ and $<_2;$ and
			\item for all $i\neq j$, the monomials $LM_{<_2}(\theta_i)$ and $LM_{<_2} (\theta_j)$ are relatively prime.
		\end{enumerate}
	\end{itemize}	
\end{theorem}
By Theorem~\ref{Gr\"obner Bases Under Composition,diff monomial ordering, main result}, the first part of this section is finished. We start now the second part which deals with the same main questions of this section and describes the studying of the same ways of generalisation but for the reduced Gr\"obner bases case. As we said before, this case is more complicated because we must solve some problems come from the remains of free variables especially with the condition: for all $g \in G$, no monomial of $g$ lies in  $\langle LT_{<}(G - \{g\}) \rangle$, which is the main condition for the reduced case. Z. Liu and M. Wang tried to solve these problems by giving two main results to describe the answer, one gives the necessary and the other gives the sufficient conditions. These two results differ by some additional condition deals with these free variables. They also gave in some corollary and remarks some forms and cases of the value of $\Theta$ which make the composition commute with our case.
The condition: $\Theta$ is a list of permuted univariate and monic polynomials, given in Theorem~\ref{Reduced Gr\"obner Bases Under Composition, main result},
has some other forms which in some cases an equivalent conditions of it. Surely, these results finished and generalised the work done by J. Guti\'errez and R. Rubio San Miguel in  \cite{GR}, but we think that it may  need some simplifications and to make these results as one result that contains all the information of the answer.
We will not give a lot of comments about the proofs and repeat ourself here since these proofs follow what is done in \cite{GR}.
We will sometimes here change the order of some results from the given order of them in the paper \cite{LZW}. For the result of necessary conditions we have:

\begin{itemize}
	\item Lemma~\ref{diff rgb give to} is just a generalisation of Lemma~\ref{Red set give To} for different polynomial rings and monomial orderings,  where the commutativity  implies the compatibility with the monomial orderings. The two proofs are similar in the method and differ in some few calculations.
	
	\item In Lemma~\ref{diff rgb give rp} and Lemma~\ref{diff rgb give not divide}, the condition: $\Theta$ is a list of permuted univariate and monic polynomials, will be replaced by another two suitable conditions to our case. Later, we will see the relation between these three conditions in some special cases. Nothing  new in the proofs.
\end{itemize}

\begin{lemma}[\cite{LZW}] \label{diff rgb give to}
	Fix a monomial ordering $<_1 $ on the polynomial ring $K[x_1, \dots, x_n]$ and a monomial ordering $<_2$ on the polynomial ring $K[y_1, \dots, y_m]$. Let $\Theta = (\theta_1, \dots, \theta_n)$ be a list of $n$ non-zero polynomials in $K[y_1, \dots, y_m]$. Let
	\begin{itemize}
		\item (A) Composition by $\Theta $ commutes with the computation of reduced Gr\"obner bases with respect to $<_1$ and $<_2$;
		\item (B) Composition by $\Theta $ is compatible with the monomial orderings $<_1$ and $<_2$.
	\end{itemize}
	Then (A)$\Rightarrow$ (B).	
\end{lemma}

\begin{lemma}[\cite{LZW}] \label{diff rgb give rp}
	Fix a monomial ordering $<_1 $ on the polynomial ring $K[x_1, \dots, x_n]$ and a monomial ordering $<_2$ on the polynomial ring $K[y_1, \dots, y_m]$. Let $\Theta = (\theta_1, \dots, \theta_n)$ be a list of $n$ non-zero polynomials in $K[y_1, \dots, y_m]$. Let
	\begin{itemize}
		\item (A) Composition by $\Theta $ commutes with the computation of reduced Gr\"obner bases with respect to $<_1$ and $<_2$;
		\item (B) For all $i\neq j$, the monomials $LM_{<_2}(\theta_i)$ and $LM_{<_2} (\theta_j)$ are relatively prime.
	\end{itemize}
	Then (A)$\Rightarrow$ (B).	
\end{lemma}

\begin{lemma}[\cite{LZW}] \label{diff rgb give not divide}
	Fix a monomial ordering $<_1 $ on the polynomial ring $K[x_1, \dots, x_n]$ and a monomial ordering $<_2$ on the polynomial ring $K[y_1, \dots, y_m]$. Let $\Theta = (\theta_1, \dots, \theta_n)$ be a list of $n$ non-zero polynomials in $K[y_1, \dots, y_m]$. Let
	\begin{itemize}
		\item (A) Composition by $\Theta $ commutes with the computation of reduced Gr\"obner bases with respect to $<_1$ and $<_2$;
		\item (B) For all $i\neq j$  and for all $\lambda_i \geq 0$, the monomial $LM_{<_2}(\theta_j)\nmid$ any term of $ \prod _{i\neq j} \theta_i^{\lambda_i}$, and $LC_{<_2}(\Theta)= (1, \dots,1)$.
	\end{itemize}
	Then (A)$\Rightarrow$ (B).	
\end{lemma}
Two remarks explain and complete what we said before about the relations between the used conditions in the previous lemmas. The first is for a special field which is the real number field. The reason for this special result is that the product of polynomials with positive coefficients from $K=\mathbb{R}$ will save all possible monomials obtained by multiplying any monomial from the first polynomial by any monomial from the second one. If some coefficients are negative or we work over other field, then this case may  not happen. The second one is for the case $m=n$.
\begin{remark}[\cite{LZW}] \label{diff rem 1}
	If $K=\mathbb{R}$, the real number field, and every coefficient of $\theta_i$ is positive, then $<_1$ is equal to $<_2$, then the following are equivalent
	\begin{enumerate}
		\item For all $i\neq j$  and for all $\lambda_i \geq 0$, the monomial $LM_{<_2}(\theta_j)\nmid$ any term of $ \prod _{i\neq j}\theta_i^{\lambda_i}$.
		\item For all $i\neq j$  and for all $\lambda_i \geq 0$, no term of $ \prod _{i\neq j}\theta_i^{\lambda_i}$ contains a variable
		appearing in $LM_{<_2}(\theta_j)$.
	\end{enumerate}
\end{remark}

\begin{remark}[\cite{LZW}] \label{rem m=n}
	If $m=n$, it is easy to see that conditions (1) and (2) of Remark~\ref{diff rem 1} are equivalent to the fact that every $\theta_i$ is a monic polynomial in one variable, that is if we add the condition: for all $i\neq j$, the monomials $LM_{<_2}(\theta_i)$ and $LM_{<_2} (\theta_j)$ are relatively prime, then we get that  $\Theta$ is a list of permuted univariate and monic polynomials.
\end{remark}

The necessity part will be given now, and this part is just the three conditions given in the previous three lemmas but in one result. Easily, it can be proved that these conditions are a sufficient condition for composition by $\Theta$ to commute with reduced Gr\"obner basis computation if and only if $m = n$ (see \cite{LZW}). For this proof, Remark~ \ref{rem m=n} will be used.
\begin{theorem}[\cite{LZW}] \label{reduced Gr\"obner Bases Under Composition,diff monomial ordering, main result}
	Fix a monomial ordering $<_1 $ on the polynomial ring $K[x_1, \dots, x_n]$ and a monomial ordering $<_2$ on the polynomial ring $K[y_1, \dots, y_m]$. Let $\Theta = (\theta_1, \dots, \theta_n)$ be a list of $n$ non-zero polynomials in $K[y_1, \dots, y_m]$. Let
	\begin{itemize}
		\item (A) Composition by $\Theta $ commutes with the computation of reduced Gr\"obner bases with respect to $<_1$ and $<_2$;
		\item (B) We have that
		\begin{enumerate}
			\item  composition by $\Theta $ is compatible with the monomial orderings $<_1$ and $<_2;$
			\item  for all $i\neq j$, the monomials $LM_{<_2}(\theta_i)$ and $LM_{<_2} (\theta_j)$ are relatively prime; and
			\item  for all $i\neq j$  and for all $\lambda_i \geq 0$, the monomial $LM_{<_2}(\theta_j)\nmid$any term of $ \prod _{i\neq j} \theta_i^{\lambda_i}$, and $LC_{<_2}(\Theta)= (1, \dots,1)$.
		\end{enumerate}
	\end{itemize}
	Then (A)$\Rightarrow$ (B).	
\end{theorem}
Since the three conditions of Theorem~\ref{reduced Gr\"obner Bases Under Composition,diff monomial ordering, main result} are not enough to map a reduced Gr\"obner basis $G$ with respect to $<_1$ to a reduced Gr\"obner basis $G\circ \Theta$ with respect to $<_2$, we need to put another conditions to make this mapping hold. But using Theorem~\ref{Gr\"obner Bases Under Composition,diff monomial ordering, main result}, with the first two conditions, we ensure that $G\circ \Theta$ is a Gr\"obner basis with respect to $<_2$ and thus we need only some conditions preserve the conditions of reduced only.
The third condition will be enough as said before for the case $m = n$, so the needed new conditions must solve the problems of free variables in the case $n \leq m$.
This condition is: if some monomial t of $\theta_i$ contains some variable $y'$ appearing in $LM_{<_2}(\theta_i)$, then $\deg_{y'}t\leq \deg_{y'}LM_{<_2}(\theta_i)$.
In other words, any variable used in the monomial $LM_{<_2}(\theta_i)$ must have a degree in $LM_{<_2}(\theta_i)$ bigger than or equal to any degree of it in any other term in $\theta_i$.
This condition will be joined later with the third condition in Corollary~\ref{reduced Gr\"obner Bases Under Composition,diff monomial ordering, coro}
to give us one condition which show us some form of $\Theta $ to be used later for giving some examples. Although, the second condition is written in this corollary,
we think that we can ignore it since the new joined condition imply it.
\begin{theorem}[\cite{LZW}] \label{reduced Gr\"obner Bases Under Composition,diff monomial ordering, main result2}
	Fix a monomial ordering $<_1 $ on the polynomial ring $K[x_1, \dots, x_n]$ and a monomial ordering $<_2$ on the polynomial ring $K[y_1, \dots, y_m]$. Let $\Theta = (\theta_1, \dots, \theta_n)$ be a list of $n$ non-zero polynomials in $K[y_1, \dots, y_m]$. Let
	\begin{itemize}
		\item (A) Composition by $\Theta $ commutes with the computation of reduced Gr\"obner bases with respect to $<_1$ and $<_2$;
		\item (B) We have that
		\begin{enumerate}
			\item  composition by $\Theta $ is compatible with the monomial orderings $<_1$ and $<_2;$
			\item  for all $i\neq j$, the monomials $LM_{<_2}(\theta_i)$ and $LM_{<_2} (\theta_j)$ are relatively prime;
			\item  for all $i\neq j$  and for all $\lambda_i \geq 0$, the monomial $LM_{<_2}(\theta_j)\nmid$any term of $ \prod _{i\neq j} \theta_i^{\lambda_i}$, and $LC_{<_2}(\Theta)= (1, \dots,1);$ and
			\item  If some monomial t of $\theta_i$ contains some variable $y'$ appearing in $LM_{<_2}(\theta_i)$, then $\deg_{y'}t\leq \deg_{y'}LM_{<_2}(\theta_i)$.
		\end{enumerate}
	\end{itemize}
	Then (B)$\Rightarrow$ (A).	
\end{theorem}
We can make use of the following corollary to construct many examples so that composition by $\Theta $ commutes with the computation of reduced Gr\"obner bases with respect to the monomial orderings $<_1$ and $<_2$.
\begin{corollary}[\cite{LZW}] \label{reduced Gr\"obner Bases Under Composition,diff monomial ordering, coro}
	Fix a monomial ordering $<_1 $ on the polynomial ring $K[x_1, \dots, x_n]$ and a monomial ordering $<_2$ on the polynomial ring $K[y_1, \dots, y_m]$. Let $\Theta = (\theta_1, \dots, \theta_n)$ be a list of $n$ non-zero polynomials in $K[y_1, \dots, y_m]$. Let
	\begin{itemize}
		\item (A) Composition by $\Theta $ commutes with the computation of reduced Gr\"obner bases with respect to $<_1$ and $<_2$;
		\item (B) We have that
		\begin{enumerate}
			\item  composition by $\Theta $ is compatible with the monomial orderings $<_1$ and $<_2;$
			\item  for all $i\neq j$, the monomials $LM_{<_2}(\theta_i)$ and $LM_{<_2} (\theta_j)$ are relatively prime;
			\item $\text{the list} \ LM_{<_2}(\Theta)$
			
			$= (y^{\lambda_1}_{\pi_1}\dots y^{\lambda_{l_1}}_{\pi_{l_1}}, y^{\lambda_{l_1+1}}_{\pi_{l_1+1}} \dots y^{\lambda_{l_1+l_2}}_{\pi_{l_1+l_2}},\dots,y^{\lambda_{l_1+\dots l_{n-1}+1}}_{\pi_{l_1+\dots l_{n-1}+1}}\dots y^{\lambda_{l_1+\dots l_{n}}}_{\pi_{l_1+\dots l_{n}}})$,
			
			for some $\pi \in S_m$ and $\lambda_i>0$, such that
			$$\theta_1 \in K[y_{\pi_1},\dots,y_{\pi_{l_1}}],$$$$ \theta_2 \in K[y_{\pi_{l_1+1}},\dots,y_{\pi_{l_1+l_2}}],\dots, \theta_n \in K[y_{\pi_{l_1+\dots l_{n-1}+1}},\dots,y_{\pi_{l_1+\dots l_{n}}}],$$
			and for any term
			$$y^{b_1}_{\pi_1}\dots y^{b_{l_1}}_{\pi_{l_1}} $$ of $\theta_1$, $$b_1\leq\lambda_1,\dots,b_{l_1}\leq\lambda_{l_1},$$
			$$ \vdots$$
			and for any term
			$$y^{c_1}_{\pi_{l_1+\dots +l_{n-1}+1}}\dots y^{c_{l_n}}_{\pi_{l_1+\dots l_{n}}} $$
			of $\theta_n$, $$c_1\leq\lambda_{l_1+\dots l_{n-1}+1},\dots,c_{l_n}\leq\lambda_{l_1+\dots l_{n}}$$.
		\end{enumerate}
	\end{itemize}
	Then (B)$\Rightarrow$ (A).	
\end{corollary}
The last remark for the case $m=n$. Note that we ignore some other remarks because we put them in some previous comments above. After the following remark, we will end this section with some examples taken from \cite{LZW}. Some of them explain why we need to divided the answer into two theorems.
\begin{remark}[\cite{LZW}]
	If $m=n$, it is easy to show that Theorem~\ref{Reduced Gr\"obner Bases Under Composition, main result} can be obtained as a corollary of Theorems~\ref{reduced Gr\"obner Bases Under Composition,diff monomial ordering, main result} and \ref{reduced Gr\"obner Bases Under Composition,diff monomial ordering, main result2} by making $<_1$ equal to $<_2$. Note that condition (4) of Theorem~\ref{reduced Gr\"obner Bases Under Composition,diff monomial ordering, main result2} is a necessary condition for composition by $\Theta$ to be commute with the computation of reduced Gr\"obner bases if and only if $m=n$.
\end{remark}

\begin{example}  \hfill
	\begin{enumerate}[leftmargin=20pt]
		\item If $m=n$ and $<_1 = <_2 $, then Examples~\ref{ex1,1}, \ref{ex2,1} and \ref{ex3,1} given at the end of Section~\ref{Gr\"obner Bases Under Composition2} can be considered as  examples of  commutative compositions that satisfy the conditions of Theorem~\ref{Gr\"obner Bases Under Composition,diff monomial ordering, main result}.
		
		\item If $m=n$ and $<_1 = (<_2)_\Theta $, then Examples~\ref{ex1,2}--\ref{ex4,2} given at the end of Section~\ref{Gr\"obner Bases Under Composition3} can be considered as  examples of  commutative compositions that satisfy the conditions of Theorem~\ref{Gr\"obner Bases Under Composition,diff monomial ordering, main result}.
		
		\item If $m=n$ and $<_1 = <_2$, then Examples~\ref{ex1,3}--\ref{ex5,3} given at the end of Section~\ref{Reduced Gr\"obner Bases Under Composition} can be considered as  examples of  commutative compositions that satisfy the conditions of Theorem~\ref{reduced Gr\"obner Bases Under Composition,diff monomial ordering, main result}.
	\end{enumerate}
\end{example}

\begin{example}[\cite{LZW}]
	Let  $<_1$ be the graded lexicographic ordering on $K[x_1,x_2]$ with $x_2 < x_1$ and $<_2$ be the graded lexicographic ordering on $K[y_1,y_2,y_3]$ with $y_3 < y_2 < y_1$. Let $$\Theta=( y_1^2 + y_2^2 +y_3, y_2y_3 +y_3),$$ then $$LM_{<_2}(\Theta)= (y_1^2,y_2y_3).$$ It  is  easy  to  verify that the composition by $\Theta$ is  compatible with  $<_1$ and $<_2$, $LM_{<_2}(\theta_1)$ and $LM_{<_2}(\theta_2)$ are  relatively prime and condition (4)  of Theorem~\ref{reduced Gr\"obner Bases Under Composition,diff monomial ordering, main result2} is  true. But $$ (\theta_1)^2 = y_1^4+ 2y_1^2y_2^2 + 2y_1^2y_3 + 2y_2^2y_3+ y_2^4 +y_3^2,$$ $LM_{<_2}(\theta_2) = y_2y_3 \mid 2y_2^2y_3$, where $2y_2^2y_3$ is a term of $ (\theta_1)^2$. Hence  condition (3) of Theorem~\ref{reduced Gr\"obner Bases Under Composition,diff monomial ordering, main result2}  is  not  true and thus $ \{\theta_1^2,  \theta_2\}$  is  not a reduced Gr\"obner basis. Hence by  Theorem~\ref{Gr\"obner Bases Under Composition,diff monomial ordering, main result}  and  Theorem~\ref{reduced Gr\"obner Bases Under Composition,diff monomial ordering, main result2},  we  have  that the composition  by $\Theta$ commutes with the computation of Gr\"obner bases with respect to $<_1$ and $<_2$, but  does not commute with the computation of reduced Gr\"obner bases with respect to $<_1$ and $<_2$.
\end{example}

\begin{example}[\cite{LZW}] \label{ex3,4}
	Let  $<_1$ be the graded lexicographic ordering on $K[x_1,x_2]$ with $x_2 < x_1$ and $<_2$ be the graded lexicographic ordering on $K[y_1,y_2,y_3]$ with $y_3 < y_2 < y_1$.
 Let $$\Theta=( y_1^2 + y_2^2 +y_3, y_2y_3 +y_3^2),$$ then $$LM_{<_2}(\Theta)= (y_1^2,y_2y_3).$$
  As above, it  is  easy  to  verify that the composition by $\Theta$ is  compatible with  $<_1$ and $<_2$, $LM_{<_2}(\theta_1)$ and
  $LM_{<_2}(\theta_2)$ are  relatively prime and condition (3) of Theorem~\ref{reduced Gr\"obner Bases Under Composition,diff monomial ordering, main result2} is  true,
   but $\theta_2$  does not  satisfy condition (4), because $t = y_3^2$ is a term  of $\theta_2$ with
   $$  \deg_{y_3} (LM_{<_2}(\theta_2))= 1 < 2 =  \deg_{y_3} (t).$$
    Now,  the set $ G = \{ x_1^2-x_1x_2, x_2^2+x_1\}$ is a reduced Gr\"obner basis
    with respect to $<_1$ but $ G\circ\Theta$ is not a reduced Gr\"obner basis with respect to $<_2$ (see  \cite{LZW}). This example shows us why condition (4)
     of Theorem~\ref{reduced Gr\"obner Bases Under Composition,diff monomial ordering, main result2} is very important and can not be deleted.
\end{example}

\begin{example}[\cite{LZW}]
	Using the same details of the last two examples but with $$\Theta=( y_1^2 + y_2^2 +1, y_2y_3 +y_3),$$ then $$LM_{<_2}(\Theta)= (y_1^2,y_2y_3).$$ Then the composition by $\Theta$ satisfies  all  the  conditions of Theorem~\ref{reduced Gr\"obner Bases Under Composition,diff monomial ordering, main result2}, so it commutes with the computation of reduced Gr\"obner bases with respect to $<_1$ and $<_2$. The set $ G\circ\Theta$, where $ G = \{ x_1^2-x_1x_2, x_2^2+x_1\}$ is given  in Example~\ref{ex3,4}, will be a reduced Gr\"obner basis with respect to $<_2$ in this example (see \cite{LZW}). The same result will happen for the set $G' = G\circ\Theta$ if we consider $<_3$ to be the graded lexicographic ordering on $K[y_1,y_2,y_3,y_4]$ with $y_4 < y_3 < y_2 < y_1$ and taking
	$$\Theta'=( y_1^2 + y_2^2 +1, y_2y_3 +y_3, y_4^2),$$ and $G'\circ\Theta'$ is a reduced Gr\"obner basis with respect to $<_3$.
\end{example}

\begin{example}[\cite{LZW}]
	Let  $<_1$ be a monomial ordering on $K[x_1,\dots,x_n]$, $<_2$ be a monomial ordering on $K[y_1,\dots,y_{2m}]$ and $$\Theta=( y_1y_2, y_3 y_4,\dots, y_{2n-1}y_{2n} ).$$ Assuming that the composition by $\Theta$ is compatible with  the monomial orderings $<_1$ and $<_2$, then the  composition  by $\Theta$ commutes with the computation of Gr\"obner bases and with the computation of reduced Gr\"obner bases with respect to $<_1$ and $<_2$.
\end{example}
\section{Homogeneous Gr\"obner Bases under Composition}\label{Homogeneous Gr\"obner Bases Under Composition}
After the work done by Hong in his first paper in 1995 and giving some answers for some open questions given be hem in his second paper in 1996 and in the work of J. Guti\'errez and R. Rubio San Miguel in their paper in 1998, and after the generalizations of these results given by Z. Liu and M. Wang in 2001, some other authors started to study the same subject of the behaviour of composition over other types of Gr\"obner bases. In this section, we will study the work of J. Liu and M. Wang  \cite{LW1} published in 2006. Note that M. Wang  is one of the authors of the paper \cite{LZW} studied in Section~\ref{Gr\"obner Bases Under Composition: Different Polynomial Rings and Monomial Orderings} but Z. Liu and J. Liu  are different authors.
In their paper, J. Liu  and M. Wang studied the question: when does homogeneous polynomial composition commute with the computation of homogeneous Gr\"obner bases with respect to a monomial ordering $<$ under the same monomial ordering? That is, which conditions we must put on the list $\Theta$ of homogeneous polynomials with the same degree to make every homogeneous Gr\"obner basis $G$ with respect to some monomial ordering $<$ for some homogeneous ideal $I$ compose using $\Theta$ to a homogeneous Gr\"obner basis $G\circ\Theta$ with respect to the same $<$ for the homogeneous ideal $\langle I\circ\Theta \rangle$. They gave a complete answer: for every homogeneous Gr\"obner basis $G$ with respect to $<$, $G \circ\Theta$ is a homogeneous Gr\"obner basis with respect to the same monomial ordering if and only if the composition by $\Theta$ is homogeneously compatible with the monomial ordering $<$ and $\Theta$  is permuted powering. All polynomials in this section will be assumed to be always homogeneous. Also $\Theta$ will be assumed to be a list of $n$ homogeneous polynomials with the same degree in the subsequent discussions. This assumption is necessary to make the composition of a homogeneous polynomial $f$ with $\Theta$  a homogeneous polynomial $f\circ \Theta$.

We will start this study by giving the following needed three definitions:
\begin{itemize}
	\item In the first one, Definition~\ref{homogeneous Gr\"obner basis  for {<}} it will be given the meaning of being a `homogeneous Gr\"obner basis' with respect to some monomial ordering $<$ for some ideal $I$. Note that in this case, $I$ must be a homogeneous ideal. Later in Lemma~\ref{Buchberger's Criterion1 for HGB}, it will be given a method of how we can test some set $G$ and determine if it is a homogeneous Gr\"obner basis or not. This test is just a modifying of Theorem~\ref{Buchberger's Criterion1} which used in the general case.
	\item In Definition~\ref{Commutativity with homo Composition}, which is similar to Definition~\ref{Commutativity with Composition}, it will be given a mathematical condition for the studied problem in this section. In other words, what is the meaning of the commutativity of composition by $\Theta$ with the computation of homogeneous Gr\"obner bases with respect to the same monomial orderings. Remember that as it was done before, its other equivalent condition will be used in the proofs.
	\item Definition~\ref{Homogeneously Compatibility with monomial ordering} is a modifying of Definition~\ref{Compatibility with monomial ordering}  by making the degrees of the two comparable monomials are the same to get a suitable condition for the homogeneous case. The reader will notice later that this is the only difference (out of the proofs) between the main result in this section and Hong's result giving in Theorem~\ref{Gr\"obner Bases Under Composition, main result}.
\end{itemize}
\begin{definition} [\cite{LW1}] \label{homogeneous Gr\"obner basis  for {<}}
	Fix a monomial ordering $<$ on the polynomial ring $K[x_1, \dots, x_n]$. A Gr\"obner basis $G$ with respect to $<$ is said to be a homogeneous Gr\"obner basis with respect to $<$ denoted by $H. G. B. _{<} (G)$ if every element of $G$ is a homogeneous polynomial. We will use the notation $H. G. B. _{<} (G, I)$ to say that $G$ is a homogeneous Gr\"obner basis  with respect to $<$ for the ideal $I =\langle G \rangle$.
\end{definition}
\begin{definition}[\cite{LW1}] \label{Commutativity with homo Composition}
	(Commutativity of Homogeneous Gr\"obner Bases under Composition) Fix a monomial ordering $<$ on the polynomial ring $K[x_1, \dots, x_n]$. We say that the composition by $\Theta $ commutes with the computation of homogeneous Gr\"obner bases with respect to $<$ if the following formula is true for $\Theta $:
$$ \ \forall I \ \forall G \ [H. G. B._< (G, I) \Rightarrow H. G. B._< (G\circ\Theta, \langle I\circ\Theta \rangle)].$$
\end{definition}
\begin{definition}[\cite{LW1}] \label{Homogeneously Compatibility with monomial ordering}
	(Homogeneously Compatibility with monomial ordering) Fix a monomial ordering $<$ on the polynomial ring $K[x_1, \dots, x_n]$. We say that the composition by $\Theta $ is homogeneously compatible with the monomial ordering $<$  if for all monomials $p$ and $q$, the following formula is true for $\Theta$:
	\[ \ \forall p \ \forall q  \ [p < q, \deg(p)= \deg(q)~\Rightarrow~ p\circ (LM_{<}(\Theta)) < q\circ (LM_{<}(\Theta))].\]
\end{definition}
The test which will be used in the proof of sufficiency part of the main theorem of this section to determine weather a given set of polynomials form a homogeneous Gr\"obner basis or not is given in the following lemma. Note that the only difference between this lemma and Theorem~\ref{Buchberger's Criterion1} is only in the word `homogeneous'.
\begin{lemma}[\cite{LW1}] \label{Buchberger's Criterion1 for HGB}
	Let $I$ be a homogeneous ideal in the polynomial ring $K[x_1, \dots, x_n]$ and fix a monomial ordering $<$. Then a basis $G = \{g_1, \dots, g_t\}$ of $I$ is a homogeneous Gr\"obner basis with respect to $<$ of $I$ if and only if for all $g_i, g_j \in G, i \neq j$, there exist homogeneous polynomials $h_1, \dots, h_t$ such that
	\begin{enumerate}
		\item  $S_{<}(g_i, g_j) = h_1g_1+ \dots +h_tg_t$.
		\item  For every $k$, either $h_k = 0$ or
		
		\hspace{3.5cm} $LT_{<}(h_k) \cdot LT_{<}(g_k) < \lcm (LT_{<}(g_i), LT_{<}(g_j))$.
	\end{enumerate}
\end{lemma}
To reach the main theorem, we must do it cross out some helpful lemmas to be used in the proof.
The author of the studied paper in this section ignored many such lemmas depending on the existing of these results in the work of Hong in \cite{Ho1}.
They only put the most important lemmas which used directly in the proof of their main result. One of these lemmas is Lemma~\ref{hTO give 1},
which is in fact a copy of Lemma~\ref{TO give 1} for the homogeneous case. It allows us to compute the leading term and monomial with respect to some monomial ordering $<$
for the homogeneous polynomial $f\circ\Theta$ for any homogeneous polynomial $f$ if the composition by $\Theta $ is homogeneously compatible with $<$.
This will be used in the proof of sufficiency part.
\begin{lemma} [\cite{LW1}] \label{hTO give 1}
	Fix a monomial ordering $<$ on the polynomial ring $K[x_1, \dots, x_n]$, and let $\Theta$ be a list of non-zero polynomials in $K[x_1, \dots, x_n]$. Let
	\begin{itemize}
		\item (A) Composition by $\Theta $ is homogeneously compatible with the monomial ordering $<$;
		\item (B) For every homogeneous polynomial $f \in K[x_1, \dots, x_n]$, we have
		\begin{enumerate}
			\item $LT_{<}(f\circ\Theta) = LT_{<}(f)\circ (LT_{<}(\Theta))$.
			\item $LM_{<}(f\circ\Theta) = LM_{<}(f)\circ (LM_{<}(\Theta))$.
		\end{enumerate}
	\end{itemize}
	Then (A) $\Rightarrow$ (B).	
\end{lemma}
After Lemma~\ref{hTO give 1}, J. Liu and M. Wang put the following result as Lemma 3.2.  of their paper \cite[page 671]{LW1}.
\begin{lemm}[\!\!{\cite[Lemma 3.2.]{LW1}}]
	Fix a monomial ordering $<$ on the polynomial ring $K[x_1, \dots, x_n]$, and let $\Theta$ be a list of non-zero polynomials in $K[x_1, \dots, x_n]$. Let
	\begin{itemize}
		\item (A) The list $LM_{<}(\Theta)$ is permuted powering;
		\item (B) Composition by $\Theta $ is homogeneously compatible with the monomial ordering $<$.
	\end{itemize}
	Then (A) $\Rightarrow$ (B).
\end{lemm}

This lemma is not true for two reasons. The first is because there are many counterexamples for it. One of them is if we take $p=x$, $q=y$, $\Theta = (y,x)$ and any monomial ordering $<$ with $x < y$.
Clearly, the list $LM_{<}(\Theta)$ is permuted powering and $p < q$. But we have that $$p\circ (LM_{<}(\Theta))=y > x= q\circ (LM_{<}(\Theta)),$$
and thus the composition by $\Theta $ is not homogeneously compatible with the monomial ordering $<$. The second reason is what they said in the proof of it: ``This is just  \cite[Lemma 4.3]{Ho1}''.

Note that   \cite[Lemma 4.3]{Ho1} is Lemma~\ref{ND give 2}  which set that if the composition by $\Theta $ is compatible with the non-divisibility with respect to $<$ (equivalent to that the list $LM_{<}(\Theta)$ is permuted powering), the condition
$$ \ \forall p \ \forall q \ [\lcm (p\circ (LM_{<}(\Theta)),q\circ (LM_{<}(\Theta)))  = \lcm (p,q) \circ (LM_{<}(\Theta))]$$ is true. This condition is what we need in the calculations of the proof of the main result here as what Hong did in \cite{Ho1} for his main result.
We think that they replaced this condition by the condition (composition by $\Theta $ is homogeneously compatible with the monomial ordering $<$) as a writing mistake.
We searched if there is any other copy or modification of the paper \cite{LW1} in any place or journal and every copy we found have the same mistake.
Now by Lemma~\ref{ND give 2} and Lemma~\ref{hTO give 1}, the proof of sufficiency part in this section will be completed using the same proof of Theorem~\ref{Gr\"obner Bases Under Composition, suff}
invented by Hong.

To get the second direction, J. Liu and M. Wang used first the following copy of Lemma~\ref{not comp give TO} for homogeneous case without proving it or setting any previous lemmas for it referring to the work of Hong. This finished one half of the necessity part.
\begin{lemma}[\cite{LW1}] \label{HGB is TO}
	Fix a monomial ordering $<$ on the polynomial ring $K[x_1, \dots, x_n]$. Let $\Theta$ be a list of non-zero polynomials in $K[x_1, \dots, x_n]$. Let
	\begin{itemize}
		\item (A) Composition by $\Theta $ commutes with the computation of homogeneous Gr\"obner bases with respect to $<$;
		\item (B) Composition by $\Theta $ is homogeneously compatible with the monomial ordering $<$.
	\end{itemize}
	Then (A) $\Rightarrow$ (B).	
\end{lemma}
The differences in the method of proofs between them and Hong will be cleared in the second half of the necessity part of their main results even that this second half is the same.
This is because we work here with homogeneous polynomials and the set which need to be used here must consist of homogeneous polynomials as we will see in Lemma~\ref{not hgb}.
Remember that, if the composition by $\Theta$ commutes with the computation of general Gr\"obner bases, then every homogeneous Gr\"obner basis will be composed to a Gr\"obner basis which need not be homogeneous. In other words, the type of commutativity studied here is not related to the general one and there is no relation between them. Note that the second half here will be using the condition (the list $LM_{<}(\Theta)$ is permuted powering) instead of the condition (the composition by $\Theta $ is compatible with the non-divisibility) as Hong done. To reach this half and in order to simplify notation, using Lemma~\ref{HGB is TO} and without loss of generality, we may assume that $$x_1 > x_2 > \dots > x_n.$$
Hence we have $$LM_{<}(\theta_1)> LM_{<}(\theta_2)> \dots > LM_{<} (\theta_n).$$
Now:
\begin{itemize}
	\item Lemma~\ref{gcd} is not related only to our subject and valid for the general case since its proof does not depend on homogeneous polynomials. It gives us a method to test weather a set of two non-zero polynomials is a Gr\"obner basis with respect to some $<$ or not using the leading monomial of their greatest common divisor. Here, it will be used only in the proof of Lemma~\ref{not hgb}
	
	\item Lemma~\ref{not hgb} will be used to prove Lemma~\ref{HGB is rp} and plays the same role played by Lemma~\ref{not gb2} in the proof of Lemma~\ref{not comp give rp}. We will use a modified copy of this lemma in our work later in the next chapter.
	
	\item In Lemma~\ref{HGB is rp}, the result given in Lemma~\ref{not hgb} will be used to prove that if the composition by $\Theta $ commutes with the computation of homogeneous Gr\"obner bases with respect to $<$, then the list $LM_{<}(\Theta)$ consists of a pair-wise relatively prime monomials.
	This is that last step before the second half.
	
	\item In Lemma~\ref{HGB is ppp}, we finish to the needed second half, that is the list $LM_{<}(\Theta)$ is permuted powering if the composition by $\Theta $ commutes with the computation of homogeneous Gr\"obner bases with respect to $<$. The fact that $LM_{<}(\theta_i)\neq 1$ for any $i$ is clearly true since from the first we assume that $\Theta$ is a list of non-constant homogeneous polynomials with the same degree and that $$LM_{<}(\theta_1)> LM_{<}(\theta_2)> \dots > LM_{<} (\theta_n)$$ assumed above. By Lemma~\ref{HGB is rp}, $LM_{<}(\theta_i)$ and $LM_{<}(\theta_j)$ are relatively prime if $i \neq j$, and thus every $LM_{<}(\theta_i)$ can only involve one variable, and for $i \neq j$, $LM_{<}(\theta_i)$ and $LM_{<}(\theta_j)$ involve different variables. Thus we conclude that $LM_{<}(\Theta)$ is permuted powering.
	
\end{itemize}

\begin{lemma}[\cite{LW1}] \label{gcd}
	Fix a monomial ordering $<$ on $K[x_1, \dots, x_n]$. Let $f$ and $g$ be  two non-zero polynomials in $K[x_1, \dots, x_n]$, and  $d = \gcd(f,g)$. Then $ \{f,g\}$ is a Gr\"obner basis with respect to $<$ if and only if  $LM_{<}\Big(\dfrac{f}{d}\Big)$ and $LM_{<}\Big(\dfrac{g}{d}\Big)$ are relatively prime, i.e. $\gcd (LM_{<}(f), LM_{<}(g))=  LM_{<}(d)$.
\end{lemma}

\begin{lemma}[\cite{LW1}] \label{not hgb}
	Fix a monomial ordering $<$ on $K[x_1, \dots, x_n]$. Let $f$ and $g$ be  two non-zero polynomials in $K[x_1, \dots, x_n]$. Assume that $LM_{<}(f)$ and $LM_{<}(g)$  are not relatively prime with $LM_{<}(g) < LM_{<}(f)$. Then we have
	\begin{enumerate}
		\item  $ \{f,g\}$ is not a Gr\"obner basis with respect to $<$, or
		\item  $ \{f^2+g^2, fg, g^3\}$ is not a Gr\"obner basis with respect to $<$.
	\end{enumerate}
\end{lemma}
\begin{lemma}[\cite{LW1}] \label{HGB is rp}
	Fix a monomial ordering $<$ on the polynomial ring $K[x_1, \dots, x_n]$. Let $\Theta$ be a list of non-zero polynomials in $K[x_1, \dots, x_n]$. Let
	\begin{itemize}
		\item (A) Composition by $\Theta $ commutes with the computation of homogeneous Gr\"obner bases with respect to $<$;
		\item (B) The monomials $LM_{<}(\theta_1), \dots,LM_{<} (\theta_n)$ are pairwise relatively prime.
	\end{itemize}
	Then (A) $\Rightarrow$ (B).	
\end{lemma}
\begin{lemma}[\cite{LW1}] \label{HGB is ppp}
	Fix a monomial ordering $<$ on the polynomial ring $K[x_1, \dots, x_n]$. Let $\Theta$ be a list of non-zero polynomials in $K[x_1, \dots, x_n]$. Let
	\begin{itemize}
		\item (A) Composition by $\Theta $ commutes with the computation of homogeneous Gr\"obner bases with respect to $<$;
		\item (B) The list $LM_{<}(\Theta)$ is permuted powering.
	\end{itemize}
	Then (A) $\Rightarrow$ (B).	
\end{lemma}
Finishing to the main result now. We said before how we can get the sufficiency part. Joining Lemma~\ref{HGB is TO} and Lemma~\ref{HGB is ppp}, we get the necessity part.
\begin{theorem}[\cite{LW1}] \label{Homogeneous Gr\"obner Bases Under Composition, main result}
	(Main theorem of homogeneously commutativity with composition  with respect to $<$) Fix a monomial ordering $<$ on the polynomial ring $K[x_1, \dots, x_n]$. Let $\Theta$ be a list of non-zero polynomials in $K[x_1, \dots, x_n]$. Then the following are equivalent
	\begin{itemize}
		\item (A) Composition by $\Theta $ commutes with the computation of homogeneous Gr\"obner bases with respect to $<$;
		\item (B) Composition by $\Theta $ is
		\begin{enumerate}
			\item homogeneously compatible with the monomial ordering $<$;  and
			\item the list $LM_{<}(\Theta)$ is permuted powering.
		\end{enumerate}
	\end{itemize}	
\end{theorem}
Some examples given by us since the authors of the studied paper did not give.
\begin{example}
	If $\Theta$ is a list of homogeneous polynomials of the same degree and
	the	composition by $\Theta $ commutes with the computation of Gr\"obner bases with respect to $<$, then it commutes with the computation of homogeneous Gr\"obner bases with respect to $<$. This is because the compatibility with some monomial ordering is more general of homogeneously compatibility with it. Thus, Example~\ref{ex1,1},  given at the end of Section~\ref{Gr\"obner Bases Under Composition2} can be modified to be an example of this section by making $\Theta$ be a list of homogeneous polynomials of the same degree.
\end{example}
\begin{example}
	Assume that the composition by $\Theta $ is homogeneously compatible with the lexicographic ordering, then it will be homogeneously compatible with the graded lexicographic ordering $<$. Thus if the list $LM_{<}(\Theta)$ is permuted powering, then the composition by $\Theta $ commutes with the computation of homogeneous Gr\"obner bases with respect to $<$.
\end{example}
\section{$\varGamma$-Homogeneous Gr\"obner Bases under Composition}\label{Further results on Gr\"obner Bases Under Composition}

J. Liu and M. Wang decided to generalise their
results \cite{LW1} and the results proved by Hong   \cite{Ho1} unifying these results in a common result. In 2007
they published the paper  \cite{LW2}, where they studied the question of when $\varGamma$-homogeneous Gr\"obner bases
under an arbitrary grading  $\varGamma$ remain Gr\"obner bases under composition with respect to any fixed monomial ordering and gave a complete answer and characterization for it.
In more details, let $\varGamma$ be an arbitrary grading on $K[x_1, \dots, x_n]$, what is the behaviour of $\varGamma$-homogeneous Gr\"obner bases under composition of polynomials, that is,
if $G$ is any $\varGamma$-homogeneous Gr\"obner basis with respect to some monomial ordering $<$ for some ideal $I$, what are the conditions must the list $\Theta$
have to make the set $G\circ\Theta$ a Gr\"obner basis with respect to the same $<$ for the ideal $\langle I\circ\Theta \rangle$.
They gave necessary and sufficient conditions to the above problem which is that the composition by $\Theta $ is $\varGamma$-compatible with the monomial ordering $<$
and the list $LM_{<}(\Theta)$ is permuted powering. If we choose some special values of $\varGamma$, then we get the results proved by Hong in \cite{Ho1} and by themselves in \cite{LW1}, which,
as we said above, may be regarded as a common generalization of these results. In this section, no other special conditions
will be put on $\Theta$ since we study the general case and the obtained composed set $G\circ\Theta$ needed to be only a Gr\"obner basis without any special cases.

As done before, we will start this study by giving the following needed three definitions:
\begin{itemize}
	\item In Definition~\ref{gamma weight} three definitions will be given. First, the meaning of being a `grading' on a polynomial ring and how we can compute some $\varGamma$-weight of a monomial $p$ for some fixed $\varGamma$. Using this $\varGamma$-weight, we can know when a polynomial $f$ is said to be $\varGamma$-homogeneous. Note that, there is no relation between the normal homogeneous polynomials and $\varGamma$-homogeneous polynomials in general. This because we can find some polynomial $f$ that is homogeneous but not $\varGamma$-homogeneous and vice versa. Finally, when a set of polynomials forms a $\varGamma$-homogeneous Gr\"obner basis  with respect to some monomial ordering $<$.
	
	\item In Definition~\ref{Commutativity with gamma homo Composition}, which is similar to Definition~\ref{Commutativity with Composition} and its copies, it will be given a mathematical condition for the studied problem in this section, i.e. the meaning of the commutativity of composition by $\Theta$ with the computation of $\varGamma$-homogeneous Gr\"obner bases with respect to the same monomial orderings. As before also, its equivalent condition will be used in the proofs. Notice that in this definition, we require that $G\circ\Theta$ is only a Gr\"obner basis and so we do not need a restriction on $\Theta$  as in \cite{LW1}.
	
	\item Definition~\ref{gamma Compatibility with monomial ordering} differs with Definition~\ref{Compatibility with monomial ordering} by  restricting it to be for monomials with equal $\varGamma$-weight only as what done in Definition~\ref{Homogeneously Compatibility with monomial ordering}. This makes  Definition~\ref{Compatibility with monomial ordering} the stronger one of them, and if it is true then both  Definitions~\ref{Homogeneously Compatibility with monomial ordering} and  \ref{gamma Compatibility with monomial ordering} will be true for any $\varGamma$. The only different between the results given in the three papers \cite{Ho1,LW1,LW2} (out of proofs) is only from these three definitions.
\end{itemize}
\begin{definition}[\cite{LW2}] \label{gamma weight}
	Let $\varGamma = (\gamma_1, \dots, \gamma_n)$ be an $n$-tuple of non-negative integers and fix a monomial ordering $<$ on the polynomial ring $K[x_1, \dots, x_n]$.
	\begin{enumerate}[leftmargin=17pt]
		\item The $\varGamma$-weight of a monomial $p = x_1^{u_1}\dots x_n^{u_n}$ is defined as $$\varGamma(p)= \gamma_1u_1+ \dots +\gamma_nu_n.$$ Under this sense, $\varGamma$ is said to be a grading on $K[x_1, \dots, x_n]$.
		\item We say that a polynomial $f$ is $\varGamma$-homogeneous if  every monomial in $f$ has the same $\varGamma$-weight.
		\item A Gr\"obner basis $G$ with respect to $<$ is said to be a $\varGamma$-homogeneous denoted by $\varGamma. H. G. B. _{<} (G)$ if every polynomial of $G$ is $\varGamma$-homogeneous. We will use the notation $\varGamma. H. G. B. _{<} (G, I)$ to say that $G$ is a $\varGamma$-homogeneous Gr\"obner basis  with respect to $<$ for the ideal $I =\langle G \rangle$.
	\end{enumerate}
\end{definition}
\begin{definition}[\cite{LW2}] \label{Commutativity with gamma homo Composition}
	(Commutativity of $\varGamma$-Homogeneous Gr\"obner Bases under Composition) Let $\varGamma = (\gamma_1, \dots, \gamma_n)$ be an $n$-tuple of non-negative integers. Fix a monomial ordering $<$ on the polynomial ring $K[x_1, \dots, x_n]$. We say that the composition by $\Theta $ commutes with the computation of $\varGamma$-homogeneous Gr\"obner bases with respect to $<$ if the following formula is true for $\Theta $: $$ \ \forall I \ \forall G  \ [\varGamma. H. G. B._< (G, I) \Rightarrow G. B._< (G\circ\Theta, \langle I\circ\Theta \rangle)].$$
\end{definition}

\begin{definition}[\cite{LW2}] \label{gamma Compatibility with monomial ordering}
	($\varGamma$-Compatibility with monomial ordering) Let $\varGamma = (\gamma_1, \dots, \gamma_n)$ be an $n$-tuple of non-negative integers. Fix a monomial ordering $<$ on the polynomial ring $K[x_1, \dots, x_n]$. We say that the composition by $\Theta $ is $\varGamma$-compatible with the monomial ordering $<$  if  for all monomials $p$ and $q$, the following formula is true for $\Theta $: $$ \ \forall p \ \forall q \ [p < q, \varGamma(p)= \varGamma(q)~\Rightarrow~ p\circ (LM_{<}(\Theta)) < q\circ (LM_{<}(\Theta))].$$
\end{definition}
The following remark describes why this study is a generalisation of \cite{Ho1} and \cite{LW1} using some special values of $\varGamma$.

\begin{remark}[\cite{LW2}] \label{special values} \hfill
	\begin{itemize}
		\item If $\varGamma (x_1)= \dots =\varGamma (x_n) = 0$, then Definition~\ref{Commutativity with gamma homo Composition} is the same as Definition~\ref{Commutativity with Composition} and Definition~\ref{gamma Compatibility with monomial ordering} is the same as Definition~\ref{Compatibility with monomial ordering}.
		
		\item If $\varGamma (x_1)= \dots =\varGamma (x_n) = 1$, then Definition~\ref{Commutativity with gamma homo Composition} is the same as Definition~\ref{Commutativity with homo Composition} and Definition~\ref{gamma Compatibility with monomial ordering} is the same as Definition~\ref{Homogeneously Compatibility with monomial ordering} (clearly $\Theta$ must be assumed to be a list of $n$ homogeneous polynomials of the same degree here).
	\end{itemize}	
\end{remark}
As J. Liu and M. Wang did in \cite{LW1}, they ignored the proofs of their results which  are similar to the proofs of some Hong's result given in \cite{Ho1}.
They started directly with the necessity part because the sufficiency part is similar to what was done in \cite{Ho1} and \cite{LW1}.
The first part of necessity will be given in Lemma~\ref{gammaGB is TO} without proof since it is similar to the proofs of Lemma~\ref{not comp give TO} and Lemma~\ref{HGB is TO}
where these two lemmas can be obtained from Lemma~\ref{gammaGB is TO} using the special values of $\varGamma$ given in Remark~\ref{special values}.
In this lemma, it will be given that if the composition by $\Theta $ commutes with the computation of $\varGamma$-homogeneous Gr\"obner basis with respect to $<$,
then this composition is $\varGamma$-compatible with $<$.
\begin{lemma}[\cite{LW2}] \label{gammaGB is TO}
	Let $\varGamma = (\gamma_1, \dots, \gamma_n)$ be an $n$-tuple of non-negative integers. Fix a monomial ordering $<$ on the polynomial ring $K[x_1, \dots, x_n]$. Let $\Theta$ be a list of non-zero polynomials in $K[x_1, \dots, x_n]$. Let
	\begin{itemize}
		\item (A) Composition by $\Theta $ commutes with the computation of $\varGamma$-homo\-geneous Gr\"obner bases with respect to $<$;
		\item (B) Composition by $\Theta $ is $\varGamma$-compatible with the monomial ordering $<$.
	\end{itemize}
	Then (A) $\Rightarrow$ (B).	
\end{lemma}
Besides Lemma~\ref{gcd}, which gives us that a set of two non-zero polynomials $\{f,g\}$ is a Gr\"obner basis with respect to some $<$  if and only if  $LM_{<}\Big(\dfrac{f}{d}\Big)$ and $LM_{<}\Big(\dfrac{g}{d}\Big)$ are relatively prime,
where $d = \gcd (f, g)$, the following lemma will be very important in the proofs of lemmas (Lemma~\ref{gammaGB is rp1}--Lemma~\ref{gammaGB is rp6}).
This lemma shows us that a Gr\"obner basis with respect to some $<$ will remain a Gr\"obner basis with respect to the same $<$ if we multiply each element of it by a fixed non-zero polynomial.
The converse is true also. As Lemma~\ref{gcd}, Lemma~\ref{GB iff hgb} is valid for the general case.
\begin{lemma}[\cite{LW2}] \label{GB iff hgb}
	Fix a monomial ordering $<$ on the polynomial ring $K[x_1, \dots, x_n]$. Let $G =\{g_1,\dots,g_t\} \subseteq K[x_1, \dots, x_n]$ and given a non-zero polynomial $h$ in $K[x_1, \dots, x_n]$.
	Then, $G =\{g_1,\dots,g_t\}$ is a Gr\"obner basis with respect to $<$ if and only if $hG =\{hg_1,\dots,hg_t\}$ is a Gr\"obner basis with respect to $<$.
	In other words: $$G. B._< (G) \Longleftrightarrow G. B._< (hG).$$
\end{lemma}
Before studying the second half of necessity, we need the following remark, which describes some useful used notation in the lemmas after it. This method of writing tell us why we need Lemma~\ref{gcd} and Lemma~\ref{GB iff hgb} in this section.
\begin{remark}[\cite{LW2}]
	Fix a monomial ordering $<$ on the polynomial ring $K[x_1, \dots, x_n]$. We fix the following notation for the convenience of discussions below: for $i \neq j$, let $d =\gcd(\theta_i,\theta_j)$. There are polynomials $\theta_i', \theta_j', h_i, h_j$ such that
$$\theta_i=d\theta_i',~~~~~~~~~\theta_j=d\theta_j',$$ $$ \theta_i'=LM_{<}(\theta_i') +h_i', ~~~~~~~~~\theta_j'=LM_{<}(\theta_j')+h_j',$$ where $$LM_{<}(h_i') < LM_{<}(\theta_i'),~~~~~~~~~LM_{<}(h_j') < LM_{<}(\theta_j').$$ Without loss of generality, we may assume that $LC_{<}(\theta_i')= 1$, for any $i$.
\end{remark}
The following lemma will be used in the proof of Lemma~\ref{gammaGB is rp3} below and its proof will depend on Lemma~\ref{gcd}. Note that the needed result for the second half of necessity is that for all $i\neq j$, the monomials $LM_{<}(\theta_i)$ and $LM_{<} (\theta_j)$ are relatively prime (which will be proved in the chain of Lemmas~\ref{gammaGB is rp2}--\ref{gammaGB is rp6}), not that for all $i\neq j$, the monomials $LM_{<}(\theta_i')$ and $LM_{<} (\theta_j')$ are relatively prime, given by this lemma.
\begin{lemma}[\cite{LW2}] \label{gammaGB is rp1}
	Let $\varGamma = (\gamma_1, \dots, \gamma_n)$ be an $n$-tuple of non-negative integers. Fix a monomial ordering $<$ on the polynomial ring $K[x_1, \dots, x_n]$. Let $\Theta$ be a list of non-zero polynomials in $K[x_1, \dots, x_n]$. Let
	\begin{itemize}
		\item (A) Composition by $\Theta $ commutes with the computation of $\varGamma$-homoge\-neous Gr\"obner bases with respect to $<$;
		\item (B) The monomials $LM_{<}(\theta_i')$ and $LM_{<} (\theta_j')$ are relatively prime for $i\neq j$.
	\end{itemize}
	Then (A) $\Rightarrow$ (B).	
\end{lemma}
Now, in the following chain of lemmas, we need to end with the fact that: if the composition by $\Theta $ commutes with the computation of $\varGamma$-homogeneous Gr\"obner basis with respect to $<$, then for all $i\neq j$, the monomials $LM_{<}(\theta_i)$ and $LM_{<} (\theta_j)$ are relatively prime.
This fact will be used to imply that the list $LM_{<}(\Theta)$ is permuted powering, which is the needed second half, as what done before in the previous sections.
What makes the proof of this result is over several lemmas and using several cases is that we work over arbitrary fixed $n$-tuple $\varGamma = (\gamma_1, \dots, \gamma_n)$
so that the values of $\gamma_i$ and $\gamma_j$ have three cases and each case has  different effects on the method of proof and on the finite set $G$ chosen to be a $\varGamma$-homogeneous Gr\"obner basis with respect to $<$
to be composed in the proof and get that $G \circ \Theta $ is a Gr\"obner basis with respect to $<$, then used this set to prove the relatively prime condition.
In every lemma of the following, it will be assumed that the composition by $\Theta $ commutes with the computation of $\varGamma$-homogeneous Gr\"obner bases with respect to $<$
and some special conditions on the values of $\gamma_i$ and $\gamma_j$. The final result in all of them is that the monomials $LM_{<}(\theta_i)$ and $LM_{<} (\theta_j)$ are relatively prime. Now we have:
\begin{itemize}
	\item Lemma~\ref{gammaGB is rp2} will study the case when the values of both $\gamma_i$ and $\gamma_j$ are zero. The calculations here have some common points with the work of Hong in \cite{Ho1} and this can understood since the zero value of $\varGamma$ moves us to the work of him.
	
	\item In Lemma~\ref{gammaGB is rp3} one value is zero and the other not (say here $\gamma_i \neq0$ and  $\gamma_j= 0$). The calculations here is new and we did not see it before. Lemma~\ref{gammaGB is rp1} will be used in the proof.
	
	\item Lemma~\ref{gammaGB is rp4} and Lemma~\ref{gammaGB is rp5} are like previous lemmas for Lemma~\ref{gammaGB is rp6}. In these three lemmas, the case when both of $\gamma_i$ and $\gamma_j$ are non-zero will be studied. Since the proof for this case is too long and has many cases and sup-cases, the author divided it into these three lemmas. From Lemma~\ref{gammaGB is rp4}, we get the condition (using the fixed assumption above)
	$$0< \gamma_j < \gamma_i,~x_j^{\gamma_i} < x_i^{\gamma_j}\Rightarrow LM_{<}(\theta_i')^{\gamma_j} \nmid LM_{<} (d)^{2\gamma_i-2\gamma_j}.$$
	Also from Lemma~\ref{gammaGB is rp5}, we get the condition
	$$ 0 < \gamma_j < \gamma_i,~x_i^{\gamma_j} < x_j^{\gamma_i} \Rightarrow LM_{<}(d)^{\gamma_i-\gamma_j}\cdot LM_{<}(\theta_j')^{\gamma_i} \nmid LM_{<} (\theta_i')^{2\gamma_j}.$$ These two conditions will be used to finish the case when the two non-zero values of $\gamma_i$ and $\gamma_j$ are non-equal by assuming without lost of generality that $ 0 < \gamma_j < \gamma_i$. The case when these values are equal will be finished directly in Lemma~\ref{gammaGB is rp6}. In the proofs of these three lemmas, the results given in Lemma~\ref{gcd} and Lemma~\ref{GB iff hgb} will be used many times. There are a lot of calculations during the proofs, but the main idea is what J. Liu  and M. Wang  used in \cite{LW1}, especially in Lemma~\ref{not hgb} and Lemma~\ref{HGB is rp}. They used many times sets $G$ with forms like the set $ \{f^2+g^2, fg, g^3\}$ used in these lemmas, and some times they used this set itself.
\end{itemize}

Adding the results given by these five lemmas together, we finish under the above assumption to that $LM_{<}(\theta_i)$ and $LM_{<} (\theta_j)$ are relatively prime for any $i\neq j$, and this means as shown in the previous sections that the list $LM_{<}(\Theta)$ is permuted powering if the composition by $\Theta $ commutes with the computation of $\varGamma$-homogeneous Gr\"obner bases with respect to $<$ and this will finish the second half of necessity as we said before.

\begin{lemma}[\cite{LW2}] \label{gammaGB is rp2}
	Let $\varGamma = (\gamma_1, \dots, \gamma_n)$ be an $n$-tuple of non-negative integers. Fix a monomial ordering $<$ on the polynomial ring $K[x_1, \dots, x_n]$. Let $\Theta$ be a list of non-zero polynomials in $K[x_1, \dots, x_n]$. Suppose that the composition by $\Theta $ commutes with the computation of $\varGamma$-homogeneous Gr\"obner bases with respect to $<$. If for $i\neq j, \gamma_i = \gamma_j= 0$, then the monomials $LM_{<}(\theta_i)$ and $LM_{<} (\theta_j)$ are relatively prime.
\end{lemma}

\begin{lemma}[\cite{LW2}] \label{gammaGB is rp3}
	Let $\varGamma = (\gamma_1, \dots, \gamma_n)$ be an $n$-tuple of non-negative integers. Fix a monomial ordering $<$ on the polynomial ring $K[x_1, \dots, x_n]$. Let $\Theta$ be a list of non-zero polynomials in $K[x_1, \dots, x_n]$. Suppose that the composition by $\Theta $ commutes with the computation of $\varGamma$-homogeneous Gr\"obner bases with respect to $<$. If $\gamma_i \neq0, \gamma_j= 0$, then the monomials $LM_{<}(\theta_i)$ and $LM_{<} (\theta_j)$ are relatively prime.
\end{lemma}

\begin{lemma}[\cite{LW2}] \label{gammaGB is rp4}
	Let $\varGamma = (\gamma_1, \dots, \gamma_n)$ be an $n$-tuple of non-negative integers. Fix a monomial ordering $<$ on the polynomial ring $K[x_1, \dots, x_n]$. Let $\Theta$ be a list of non-zero polynomials in $K[x_1, \dots, x_n]$. Suppose that the composition by $\Theta $ commutes with the computation of $\varGamma$-homogeneous Gr\"obner basis with respect to $<$. If $0< \gamma_j < \gamma_i$ and $x_j^{\gamma_i} < x_i^{\gamma_j}$, then  $$LM_{<}(\theta_i')^{\gamma_j} \nmid LM_{<} (d)^{2\gamma_i-2\gamma_j}.$$
\end{lemma}

\begin{lemma}[\cite{LW2}] \label{gammaGB is rp5}
	Let $\varGamma = (\gamma_1, \dots, \gamma_n)$ be an $n$-tuple of non-negative integers. Fix a monomial ordering $<$ on the polynomial ring $K[x_1, \dots, x_n]$. Let $\Theta$ be a list of non-zero polynomials in $K[x_1, \dots, x_n]$. Suppose that the composition by $\Theta $ commutes with the computation of $\varGamma$-homogeneous Gr\"obner bases with respect to $<$. If $0< \gamma_j < \gamma_i$ and $x_i^{\gamma_j} < x_j^{\gamma_i}$, then  $$ LM_{<} (d)^{\gamma_i-\gamma_j}\cdot LM_{<}(\theta_j')^{\gamma_i} \nmid LM_{<} (\theta_i')^{2\gamma_j}.$$
\end{lemma}

\begin{lemma}[\cite{LW2}] \label{gammaGB is rp6}
	Let $\varGamma = (\gamma_1, \dots, \gamma_n)$ be an $n$-tuple of non-negative integers. Fix a monomial ordering $<$ on the polynomial ring $K[x_1, \dots, x_n]$. Let $\Theta$ be a list of non-zero polynomials in $K[x_1, \dots, x_n]$. Suppose that the composition by $\Theta $ commutes with the computation of $\varGamma$-homogeneous Gr\"obner bases with respect to $<$. If $ \gamma_i\neq0$ and $ \gamma_j\neq  0$, then the monomials $LM_{<}(\theta_i)$ and $LM_{<} (\theta_j)$ are relatively prime.
\end{lemma}
Adding the results described above in Lemma~\ref{gammaGB is TO} and the chain of  Lemmas~\ref{gammaGB is rp2}--\ref{gammaGB is rp6} besides the sufficiency part which the author ignored its proof, we get the following main result of this section which is similar to the main results given in \cite{Ho1} in and \cite{LW1}.
\begin{theorem}[\cite{LW2}] \label{gamma Homogeneous Gr\"obner Bases Under Composition, main result}
	(Main theorem of $\varGamma$-homogeneously commutativity with composition  with respect to $<$ ) Let $\varGamma = (\gamma_1, \dots, \gamma_n)$ be an $n$-tuple of non-negative integers. Fix a monomial ordering $<$ on the polynomial ring $K[x_1, \dots, x_n]$. Let $\Theta$ be a list of non-zero polynomials in $K[x_1, \dots, x_n]$, then the following are equivalent
	\begin{itemize}
		\item (A) Composition by $\Theta $ commutes with the computation of $\varGamma$-homo\-geneous Gr\"obner bases with respect to $<$;
		\item (B) Composition by $\Theta $ is
		\begin{enumerate}
			\item $\varGamma$-compatible with the monomial ordering $<$;  and
			\item the list $LM_{<}(\Theta)$ is permuted powering.
		\end{enumerate}
	\end{itemize}	
\end{theorem}
\begin{remark}[\cite{LW2}]\label{sp of gamma}
	Using the special values of $\varGamma$ given in Remark~\ref{special values}, we have that:
	\begin{itemize}
		\item If $\varGamma (x_1)= \dots =\varGamma (x_n) = 0$, then Theorem~\ref{gamma Homogeneous Gr\"obner Bases Under Composition, main result} is the same as Theorem~\ref{Gr\"obner Bases Under Composition, main result}.
		
		\item If $\varGamma (x_1)= \dots =\varGamma (x_n) = 1$, then Theorem~\ref{gamma Homogeneous Gr\"obner Bases Under Composition, main result} is the same as Theorem~\ref{Homogeneous Gr\"obner Bases Under Composition, main result} (where $\Theta$ is a list of $n$ homogeneous polynomials with the same degree).
	\end{itemize}
\end{remark}
We finish with some examples.

\begin{example}
	Using Remark~\ref{sp of gamma}, Examples~\ref{ex1,1}, \ref{ex2,1} and \ref{ex3,1} given at the end of Section~\ref{Gr\"obner Bases Under Composition2} are examples of this section with $\varGamma (x_1)= \dots =\varGamma (x_n) = 0$. Same thing with the examples given at the end of Section~\ref{Homogeneous Gr\"obner Bases Under Composition} with $\varGamma (x_1)= \dots =\varGamma (x_n) = 1$.
	
\end{example}
\begin{example}
	Recall the monomial ordering $<_L$ defined in Example~\ref{ex3,3} where the choosing numbers $\alpha_1, \dots, \alpha_n$ are positive integers and let
	$\varGamma = (\alpha_1, \dots, \alpha_n)$. Then, if the composition by $\Theta $ commutes with the computation of Gr\"obner basis with respect to $<$, then it commutes with the computation of $\varGamma$-homogeneous Gr\"obner bases with respect to $<_L$.
	
\end{example}


\end{document}